\documentclass[a4paper,12pt]{article}
\usepackage[utf8]{inputenc}
\usepackage{amssymb,latexsym,amsmath}
\usepackage{amsthm}
\usepackage{tikz}
\usetikzlibrary{arrows}
\usepackage{graphicx,times} 
\usepackage{hyperref}

\newtheorem{exa}{Example}[section]
\newtheorem{Def}{Definition}[section]

\newtheorem{cor}{Corollary}[section]
\newtheorem{Pro}{Proposition}[section]
\newtheorem{conj}{Conjecture}[section]
\newtheorem{rem}{Remark}[section]
\title{Zeta functions of finite Schreier graphs and their zig zag products}
\author{Asif Shaikh, 
Hemant Bhate}
\date{}
\begin{document}
\maketitle
\begin{abstract}
We investigate the Ihara zeta functions of finite Schreier graphs $\Gamma_n$ of the Basilica group. We show that $\Gamma_{1+n}$ is $2$ sheeted unramified normal covering of $\Gamma_n, ~\forall~ n \geq 1$ with Galois group $\displaystyle \frac{\mathbb{Z}}{2\mathbb{Z}}.$ In fact, for any $n > 1, r \geq 1$ the graph $\Gamma_{n+r}$ is $2^n$ sheeted unramified, non normal covering of $\Gamma_r.$  
In order to do this we give the definition of the $generalized$ $replacement$ $product$ of Schreier graphs. We also show the corresponding results in zig zag product of Schreier graphs $\Gamma_n$ with a $4$ cycle.
\end{abstract}
\section{INTRODUCTION}
The main aim of this paper is to study the Ihara zeta functions of finite Schreier graphs $\Gamma_n$ of the Basilica group as well as Ihara zeta functions of zig zag products of these Schreier graphs with a $4$ cycle. The Ihara zeta function is a graph theoretic analogue of the Riemann zeta function $\zeta(s).$ It is defined for a connected graph $G = (V,E)$, for $t \in \mathbb{C},$ with $|t|$ sufficiently small, by 
\begin{equation}
\zeta_G(t) = \prod_{[C]~prime~cycle~in~G} (1 - t^{\nu (C)})^{-1}
\end{equation} 
where a prime $[C]$ in $G$ is an equivalence class of tailless, backtrackless primitive cycles $C$ in $G$, length of $C$ is $\nu (C).$ See \cite{Ter}. The connection of $\zeta_G(t)$ with the adjacency matrix $A$ of $G$ is given by Ihara-Bass's Theorem which says
\begin{equation}
 \zeta_G(t)^{-1} = (1-t^2)^{r-1} \det (I - At + Qt^2),
\end{equation}
where $r = |E| - |V| - 1 = $ rank of fundamental group of $G$ and $Q$ is the diagonal matrix whose $j^{th}$ diagonal entry is $Q_{jj} = d(v_j) - 1$, where $d(v_j)$ is the degree of the $j^{th}$ vertex of $G.$ The proof of the Ihara-Bass's theorem can be found in \cite{Ter}. 
We say that finite graph $\widetilde{G}$ is an unramified covering of a finite graph $G$ if there is a covering map $\pi: \widetilde{G} \rightarrow G$ which is an onto map such that for every vertex $u \in G$ and for every $v \in \pi^{-1}(u),$ the set of points adjacent to $v$ in $\widetilde{G}$ is mapped by $\pi$ one-to-one, onto the vertices in $G$ which are adjacent to $u.$ A $d$-sheeted covering is a normal or Galois covering iff there are $d$ graph automorphisms $\sigma : \widetilde{G} \rightarrow \widetilde{G}$ such that $\pi(\sigma(v)) = \pi(v),~\forall~v \in \widetilde{G}.$ These automorphisms form the Galois group $\mathbb{G} = Gal(\widetilde{G}|G).$ See \cite{Ded,Ter} for more details.\\
Our motivation comes from a question raised by A. Terras \cite{Ter} of how zeta functions behave with respect to graph products, in particular, zig zag products. We have been able to prove the following:\\
\textit{ If $n > 1,$ then $\Gamma_{n+r}\textcircled{z} C_4 $ is the non normal covering of the graph $\Gamma_{r}  \textcircled{z}C_4.$ In fact if $n = 1$ then $\Gamma_{1+r}\textcircled{z} C_4 $ is $2$ sheeted normal covering of the graph $\Gamma_{r} \textcircled{z}C_4.$ Thus $\zeta_{(\Gamma_{r} \textcircled{z}C_4)}(u)^{-1}$ divides $\zeta_{(\Gamma_{n+r} \textcircled{z}C_4)}(u)^{-1}, ~\forall n.$}\\ In order to prove the above result a new \textit{generalized replacement product $\Gamma_n \textcircled{g} \Gamma_r$} has been defined for two Schreier graphs $\Gamma_n, \Gamma_r$ which gives the resultant graph $\Gamma_n \textcircled{g} \Gamma_r$ is also the Schreier graph $\Gamma_{n+r}.$ Also,\\
\textit{If $n>1,$ then $\Gamma_{n+r}$ is non normal covering of the graph $\Gamma_r$ . In fact $\Gamma_{1+r}$ is $2$ sheeted normal covering of the graph $\Gamma_r$ and $\zeta_{\Gamma_{r}}(u)^{-1}$ divides $\zeta_{\Gamma_{n+r}}(u)^{-1}, ~\forall n.$}\\
This paper is organized as follows: In section 2, we define Schreier graphs of the Basilica group
and give several examples. This section also contains the definition of \textit{generalized replacement product} and the computations of Ihara zeta functions of these Schreier graphs using Artin $L$ functions. This section ends with some results which will be needed in the next section. Section 3 starts with the definition of zig zag product of two graphs see \cite{Rei,DAng}. Computations of Ihara zeta functions of zig zag product of Schreier graphs with a $4$ cycle are also presented. We also prove the corresponding results in zig zag product of graphs. 
\section{SCHREIER GRAPHS OF THE BASILICA GROUP}
Let $X = \{0,1\}$ be a binary alphabet. Denote $X^0$ the set consisting of the empty word and by $ X^n = \{ w = x_1x_2\cdots x_n : x_i \in X \}$ the set of words of length $n$ over the alphabet $X,$ for each $n \geq 1.$ The Basilica group $B$ acting on $X^n$ is the group generated by three state automaton. The states $a$ and $b$ of the automaton are the generators of the group. The action of $a$ and $b$ is given by 
$$ a(0w) = 0b(w), a(1w) = 1w, b(0w) = 1a(w), b(1w) = 0w, ~~\forall~~w \in X^n.$$
For each $n \geq 1,$ let $\Gamma_n$ be the $ Schreier ~~graph~~$ associated with the action of $B$ on $X^n$ with $V(\Gamma_n) = X^n$ and two vertices $v, v'$ are are connected by an edge labeled by $s$ near $v$ and by $s^{-1}$ near $v'$ if $s(v) = v'~~(i.e.~~s^{-1}(v') = v)$, where $s \in \{a^{\pm 1}, b^{\pm 1}\}.$ As the action of $B$ on $X^n$ is transitive, the graph $\Gamma_n$ is a $4$ regular connected graph on $2^n$ vertices and labels near $v$ are given by  $a^{\pm 1}, b^{\pm 1}$ for every $v \in X^n.$ We denote $s(v)$ by $v^s,$ 
where $v \in X^n, s \in \{ a^{\pm 1}, b^{\pm 1} \}.$\\
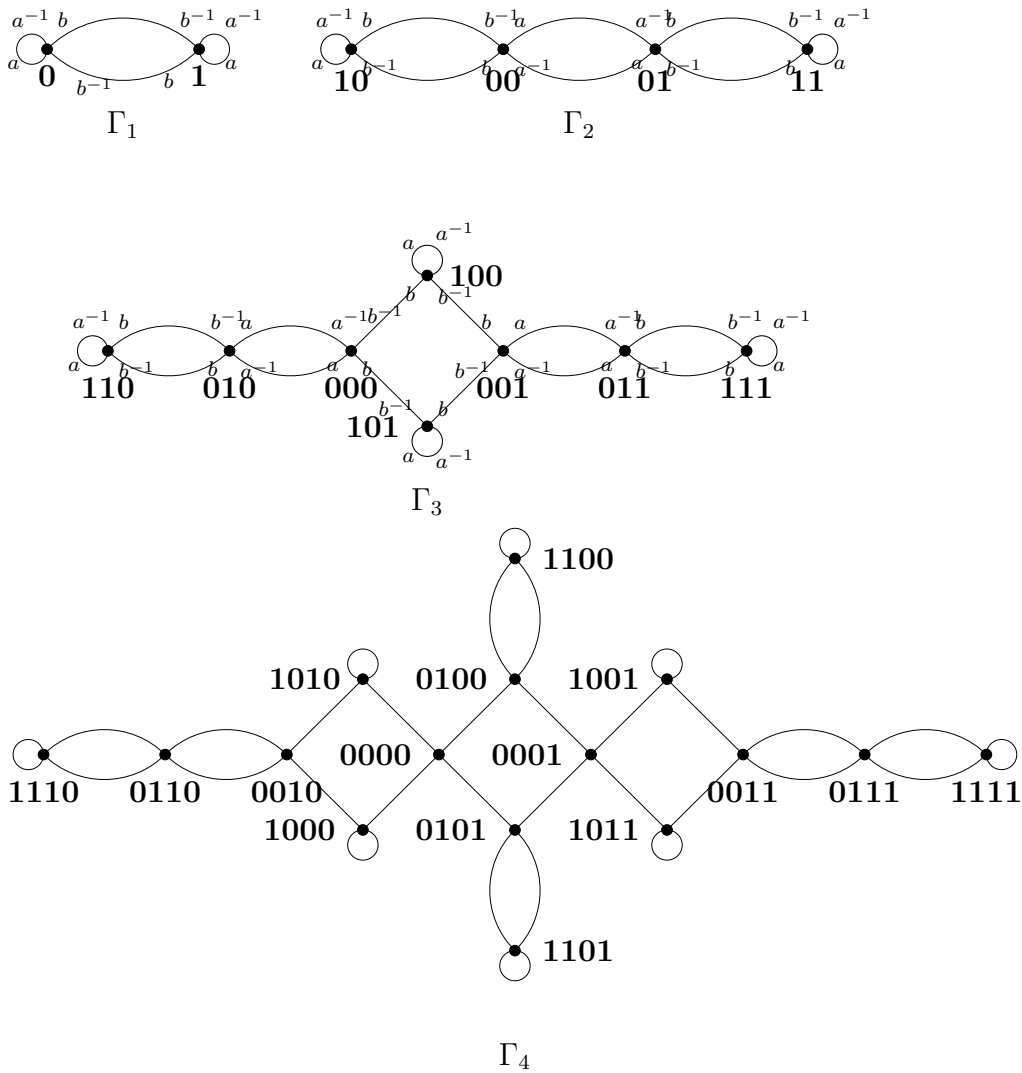
\begin{figure}[h]
\begin{tikzpicture}[scale=2]
\draw[ultra thin] (0,1) node[above right] {$^b$} to [out=45,in=135] (1,1) node[above] {$^{b^{-1}}$} ;
\draw[ultra thin] (1,1) node[below left=5pt] {$_{b}$} to [out=225,in=315] (0,1) node[below right=7pt] {$_{b^{-1}}$};
\draw[ultra thin] (1.1,1)node[above right] {$^{a^{-1}}$} circle [radius=0.1] node[below right] {$_{a}$};
\draw[ultra thin] (-0.1,1)node[above] {$^{a^{-1}}$} circle [radius=0.1] node[below left] {$_{a}$};
\draw[fill] (0,1) circle [radius=0.035] node[below=2pt]  {${\bf 0}$};
\draw[fill] (1,1) circle [radius=0.035] node[below=2pt]  {${\bf 1}$};
\draw[] (0.5,0.5) node[]{$\Gamma_1$}; 
\draw[ultra thin] (2,1) node[above right] {$^b$} to [out=45,in=135] (3,1) node[above] {$^{b^{-1}}$} ;
\draw[ultra thin] (3,1) node[below left] {$_{b}$} to [out=225,in=315] (2,1) node[below right] {$_{b^{-1}}$};
\draw[ultra thin] (3,1) node[above right] {$^a$} to [out=45,in=135] (4,1) node[above] {$^{a^{-1}}$} ;
\draw[ultra thin] (4,1) node[below left] {$_{a}$} to [out=225,in=315] (3,1) node[below right] {$_{a^{-1}}$};
\draw[ultra thin] (4,1) node[above right] {$^b$} to [out=45,in=135] (5,1) node[above] {$^{b^{-1}}$} ;
\draw[ultra thin] (5,1) node[below left] {$_{b}$} to [out=225,in=315] (4,1) node[below right] {$_{b^{-1}}$};

\draw[] (5.1,1)node[above right] {$^{a^{-1}}$} circle [radius=0.1] node[below right] {$_{a}$};
\draw[] (1.9,1)node[above] {$^{a^{-1}}$} circle [radius=0.1] node[below left] {$_{a}$};
\draw[fill] (2,1) circle [radius=0.035] node[below=4pt]  {${\bf 10}$};
\draw[fill] (3,1) circle [radius=0.035] node[below=4pt]  {${\bf 00}$};
\draw[fill] (4,1) circle [radius=0.035] node[below=4pt]  {${\bf 01}$};
\draw[fill] (5,1) circle [radius=0.035] node[below=4pt]  {${\bf 11}$};
\draw[] (3.5,0.5) node[]{$\Gamma_2$}; 
\draw[ultra thin] (0.4,-1) node[above right] {$^b$} to [out=45,in=135] (1.2,-1) node[above] {$^{b^{-1}}$} ;
\draw[ultra thin] (1.2,-1) node[below left] {$_{b}$} to [out=225,in=315] (0.4,-1) node[below right] {$_{b^{-1}}$};
\draw[ultra thin] (1.2,-1) node[above right] {$^a$} to [out=45,in=135] (2,-1) node[above] {$^{a^{-1}}$} ;
\draw[ultra thin] (2,-1) node[below left] {$_{a}$} to [out=225,in=315] (1.2,-1) node[below right] {$_{a^{-1}}$};

\draw[ultra thin] (2,-1) node[above right=2pt] {$^{b^{-1}}$} to [out=45,in=225] (2.5,-0.5) node[below left] {$^b$};         
\draw[ultra thin] (2.5,-0.5) node[below right] {$^{b^{-1}}$} to [out=-45,in=135] (3,-1) node[above left] {$^b$};
\draw[ultra thin] (3,-1) node[below left] {$_{b^{-1}}$} to [out=225,in=45] (2.5,-1.5) node[above right] {$_{b}$}; 
\draw[ultra thin] (2.5,-1.5) node[above left] {$_{b^{-1}}$} to [out=135,in=-45] (2,-1) node[below right] {$_{b}$};

\draw[ultra thin] (3,-1) node[above right] {$^a$} to [out=45,in=135] (3.8,-1) node[above] {$^{a^{-1}}$} ;
\draw[ultra thin] (3.8,-1) node[below left] {$_{a}$} to [out=225,in=315] (3,-1) node[below right] {$_{a^{-1}}$};
\draw[ultra thin] (3.8,-1) node[above right] {$^b$} to [out=45,in=135] (4.6,-1) node[above] {$^{b^{-1}}$} ;
\draw[ultra thin] (4.6,-1) node[below left] {$_{b}$} to [out=225,in=315] (3.8,-1) node[below right] {$_{b^{-1}}$};

\draw[] (4.7,-1) node[above right] {$^{a^{-1}}$} circle [radius=0.1] node[below right] {$_{a}$};
\draw[] (0.3,-1) node[above] {$^{a^{-1}}$} circle [radius=0.1] node[below left] {$_{a}$};
\draw[] (2.5,-0.4) node[above right=-1pt] {$^{a^{-1}}$} circle [radius=0.1] node[above left] {$_{a}$};
\draw[] (2.5,-1.6) node[below right=-1pt ] {$^{a^{-1}}$} circle [radius=0.1] node[below left] {$_{a}$};

\draw[fill] (.4,-1) circle [radius=0.035] node[below=6pt]  {${\bf 110}$};
\draw[fill] (1.2,-1) circle [radius=0.035] node[below=6pt]  {${\bf 010}$};
\draw[fill] (2,-1) circle [radius=0.035] node[below=6pt]  {${\bf 000}$};
\draw[fill] (2.5,-0.5) circle [radius=0.035] node[right=4pt]  {${\bf 100}$};
\draw[fill] (2.5,-1.5) circle [radius=0.035] node[left=6pt]  {${\bf 101}$};
\draw[fill] (3,-1) circle [radius=0.035] node[below=6pt]  {${\bf 001}$};
\draw[fill] (3.8,-1) circle [radius=0.035] node[below=6pt]  {${\bf 011}$};
\draw[fill] (4.6,-1) circle [radius=0.035] node[below=6pt]  {${\bf 111}$};
\draw[] (2.5,-2) node[]{$\Gamma_3$}; 
\end{tikzpicture}\\
\begin{tikzpicture}[scale=2]
\draw[ultra thin]  (0.4,-4.5)to [out=45,in=135] (1.2,-4.5) ;
\draw[ultra thin]  (1.2,-4.5)to [out=225,in=315] (0.4,-4.5) ;
\draw[ultra thin]  (1.2,-4.5)to [out=45,in=135] (2,-4.5)  ;
\draw[ultra thin]  (2,-4.5)to [out=225,in=315] (1.2,-4.5) ;

\draw[ultra thin]  (2,-4.5) to [out=45,in=225] (2.5,-4) ;         
\draw[ultra thin]  (2.5,-4) to [out=-45,in=135] (3,-4.5);
\draw[ultra thin]  (3,-4.5) to [out=225,in=45] (2.5,-5); 
\draw[ultra thin]  (2.5,-5) to [out=135,in=-45] (2,-4.5) ;

\draw[ultra thin]  (3,-4.5) to [out=45,in=225] (3.5,-4);
\draw[ultra thin]  (3.5,-4) to [out=-45,in=135] (4,-4.5);
\draw[ultra thin]  (4,-4.5) to [out=225,in=45] (3.5,-5); 
\draw[ultra thin]  (3.5,-5) to [out=135,in=-45] (3,-4.5) ;

\draw[ultra thin]  (4,-4.5) to [out=45,in=225] (4.5,-4);
\draw[ultra thin]  (4.5,-4) to [out=-45,in=135] (5,-4.5);
\draw[ultra thin]  (5,-4.5) to [out=225,in=45] (4.5,-5); 
\draw[ultra thin]  (4.5,-5) to [out=135,in=-45] (4,-4.5) ;

\draw[ultra thin]  (5,-4.5) to [out=45,in=135] (5.8,-4.5) ;
\draw[ultra thin]  (5.8,-4.5)to [out=225,in=315] (5,-4.5) ;

\draw[ultra thin] (3.5,-4) to [out=135,in=-135] (3.5,-3.2) to [out=-45,in=45] (3.5,-4);
\draw[ultra thin] (3.5,-5) to [out=-135,in=135] (3.5,-5.8) to [out=45,in=-45] (3.5,-5);

\draw[ultra thin]  (5.8,-4.5) to [out=45,in=135] (6.6,-4.5) ;
\draw[ultra thin]  (6.6,-4.5)to [out=225,in=315] (5.8,-4.5) ;

\draw[ultra thin] (6.7,-4.5) circle [radius=0.1];
\draw[ultra thin]  (0.3,-4.5) circle [radius=0.1] ;
\draw[ultra thin]  (2.5,-3.9) circle [radius=0.1] ;
\draw[ultra thin]  (2.5,-5.1) circle [radius=0.1] ;

\draw[ultra thin] (3.5,-3.1) circle [radius=0.1];
\draw[ultra thin] (3.5,-5.9) circle [radius=0.1];
\draw[ultra thin] (4.5,-3.9) circle [radius=0.1];
\draw[ultra thin] (4.5,-5.1) circle [radius=0.1];

\draw[fill] (.4,-4.5) circle [radius=0.035] node[below=6pt]  {${\bf 1110}$};
\draw[fill] (1.2,-4.5) circle [radius=0.035] node[below=6pt]  {${\bf 0110}$};
\draw[fill] (2,-4.5) circle [radius=0.035] node[below=6pt]  {${\bf 0010}$};
\draw[fill] (2.5,-4) circle [radius=0.035] node[left=4pt]  {${\bf 1010}$};
\draw[fill] (2.5,-5) circle [radius=0.035] node[left=6pt]  {${\bf 1000}$};
\draw[fill] (3,-4.5) circle [radius=0.035] node[left=6pt]  {${\bf 0000}$};
\draw[fill] (3.5,-4) circle [radius=0.035] node[left=6pt]  {${\bf 0100}$};
\draw[fill] (4,-4.5) circle [radius=0.035] node[left=6pt]  {${\bf 0001}$};
\draw[fill] (3.5,-5) circle [radius=0.035] node[left=6pt]  {${\bf 0101}$};
\draw[fill] (4.5,-4) circle [radius=0.035] node[left=6pt]  {${\bf 1001}$};
\draw[fill] (5,-4.5) circle [radius=0.035] node[below=6pt]  {${\bf 0011}$};
\draw[fill] (4.5,-5) circle [radius=0.035] node[left=6pt]  {${\bf 1011}$};
\draw[fill] (5.8,-4.5) circle [radius=0.035] node[below=6pt]  {${\bf 0111}$};
\draw[fill] (6.6,-4.5) circle [radius=0.035] node[below=6pt]  {${\bf 1111}$};
\draw[fill] (3.5,-3.2) circle [radius=0.035] node[right=6pt]  {${\bf 1100}$};
\draw[fill] (3.5,-5.8) circle [radius=0.035] node[right=6pt]  {${\bf 1101}$};
\draw[] (3.5,-6.5) node[]{$\Gamma_4$};
\end{tikzpicture}
\caption{The graph $\Gamma_1,\Gamma_2, \Gamma_3$ and $\Gamma_4$ are the Schreier graphs of the Basilica group over $X, X^2, X^3$ and $X^4$ respectively.}
\label{fig1}
\end{figure}
More information about such finite and infinite Schreier graphs of Basilica group can be found in \cite{DAn,DAng}. The complete classification (up to isomorphism) of the limiting case of infinite Schreier graphs associated with the Basilica group acting on the binary tree, in terms of the infinite binary sequence is given by D. D'Angeli, A. Donno, M. Matter and T. Nognibeda\cite{DAn}. D. D'Angeli, A. Donno and E. Sava-Huss\cite{DAng} have used Schreier graphs of Basilica groups in computations of zig zag product of graphs. The Basilica group belongs to the important class of self-similar groups and was introduced by R. Grigorchuk and A. $\dot{Z}$uk\cite{Gri}.
{\exa  The first four finite Schreier graphs of Basilica group are as in Figure 1.}\\[5pt] 
The substitution rules \cite{DAng,Grig} can be used to construct the Schreier graph $\Gamma_{n+1}$ from the graph $\Gamma_n.$ In the next sub section we give the new definition of the product of two Schreier graphs, we call it as $generalized~~replacement~~product$ of Schreier graphs. 
\subsection{Generalized replacement product}
The replacement product of two graphs is well known in literature. If $G_1$ and $G_2$ are two regular graphs with the regularity $d_1$ and $d_2$ respectively, then their replacement product $G_1 \textcircled{r} G_2$ is again a regular graph with regularity $d_2 + 1.$ The details about this product can be found in \cite{DAng,Abd}. Below is the generalization of this product, in which both $G_1$ and $G_2$ along with the resultant graph have same regularity. The motivation of generalized version is the following:  Using any two Schreier graphs of the Basilica group we should be able to produce the resultant graph to be Schreier graph of the Basilica group.\\
\\   
Let $\Gamma_n$ and $\Gamma_r$ be two Schreier graphs of Basilica group. To define their generalized replacement product, we first choose the spanning subgraph $\Gamma_r'$ of the graph $\Gamma_r$ such that $E(\Gamma_r') = E(\Gamma_r) \backslash \{e_a, e_b\}$ where $e_a$ and $e_b$ are the edges of the vertex $u_0 = 0^r \in \Gamma_r$ with labels $a$ and $b$ near $u_0$ respectively. Note that the vertex $u_0 = 0^r$ we mean the word of length $r$ and having all alphabets are $0$.
Let $a(u_0) = u_0^a$ and $b(u_0) = u_0^b$ be the vertices of $\Gamma_r'$ which are adjacent to $u_0$ in $\Gamma_r$ by the edges  $e_a$ and $e_b$ respectively.\\
Thus every vertex from the set $X^r \backslash \{u_0, u_0^a, u_0^b\}$ has degree $ d = 4$ and 
$$d(u_0) = 2, d(u_0^a) = d(u_0^b) = 3$$
$$ (4 - d(u_0)) + (4 - d(u_0^a)) + (4 - d(u_0^b)) = 2 + 1 + 1 = 4 $$
Note that in above equation sum of LHS should be the regularity $d$ of the graph. \\
\\
Notice also that in the Schreier graphs $\Gamma_n$ of the Basilica group, if $e = \{v,v'\}$ is an edge which has color say $s$ near $v$ and $s^{-1}$ near $v'$, then the $rotation~~map$ 
${\bf Rot}: X^n \times D \rightarrow X^n \times D$ is defined by 
$$ {\bf Rot}_{\Gamma_n}(v,s) = (v',s^{-1}), ~~\forall~~v,v' \in X^n,~~s,s^{-1} \in D = \{a^{\pm 1}, b^{\pm 1}\}.$$

{\Def 
The generalized replacement product $\Gamma_n \textcircled{g} \Gamma_r$ is the $4$ regular graph with vertex set $X^{n+r} = X^n  \times X^r,$ and whose edges are described by the following rotation map: $\forall ~~ v \in X^n$
\begin{eqnarray}
 {\bf Rot}_{\Gamma_n \textcircled{g} \Gamma_r}((v,u_0),a) = \left\{              
                    \begin{array}{ll}           
                   ((v^a,u_0^a),a^{-1}) & \mathrm{if~r~is~even} \\
                    ((v^b,u_0^a),a^{-1}) & \mathrm{if~r~is~odd} \\ 
                    \end{array}       
                    \right. \\
{\bf Rot}_{\Gamma_n \textcircled{g} \Gamma_r}((v,u_0),b) = \left\{              
                    \begin{array}{ll}           
                   ((v^b,u_0^b),b^{-1}) & \mathrm{if~r~is~even} \\
                   ((v^a,u_0^b),b^{-1}) & \mathrm{if~r~is~odd} \\                 
                    \end{array}       
                    \right. \\
{\bf Rot}_{\Gamma_n \textcircled{g} \Gamma_r}((v,u),s) = 
                    ((v, u^s),s^{-1}),
                    ~~\forall ~~ u \in X^r \backslash \{u_0\}~~s \in  \{ a^{\pm 1}, b^{\pm 1} \}
 \end{eqnarray} }
 We call  edges given in (3) as lifts of the edge $e_a$ and in (4) as lifts of the edge $e_b.$\\
 \\
 One can imagine that the vertex set $X^{n+r}$ of the graph $\Gamma_n \textcircled{g} \Gamma_r$ is partitioned into the sheets, which are indexed by the vertices of $\Gamma_n,$ where by definition of the $v^{th}$ sheet, for $v \in X^n,$ consists of the vertices $\{(v,u)| u \in X^r\}.$ Within this construction the idea is to put the copy of $\Gamma_r'$ around each vertex $v$ of $\Gamma_n,$ while keeping all the edges of $\Gamma_r'$ as it is and the $\Gamma_n$ can be thought of using the above given rotations in the equations (3) and (4).  The connectedness of $\Gamma_r'$ and $\Gamma_n$ guarantees about the connectedness of the graph $\Gamma_n \textcircled{g} \Gamma_r.$\\[5pt]
If an edge $e = \{v,v'\}$ in $\Gamma_n$ has colors $s$ and $s^{-1}$ near $v$ and $v'$ respectively, then we call the edge $e$ as $s$-edge. Hence we can say that in the graph $\Gamma_n$, there are two types of edges i.e. $a$-edges and $b$-edges.   
 \\
 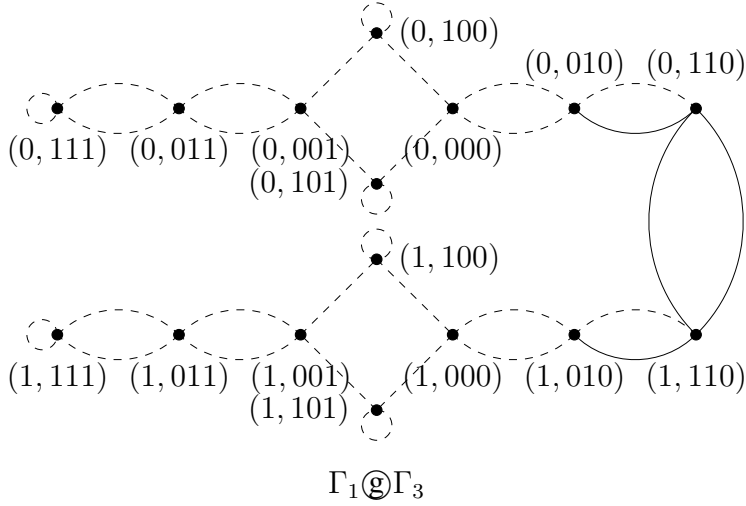
\begin{figure}[h]
\begin{tikzpicture}[scale=2]
\draw[dashed, ultra thin] (0.4,0)to [out=45,in=135] (1.2,0) ;
\draw[dashed, ultra thin] (1.2,0)to [out=225,in=315] (0.4,0) ;
\draw[dashed, ultra thin] (1.2,0)to [out=45,in=135] (2,0)  ;
\draw[dashed, ultra thin] (2,0)to [out=225,in=315] (1.2,0) ;

\draw[dashed, ultra thin] (2,0) to [out=45,in=225] (2.5,0.5) ;         
\draw[dashed, ultra thin] (2.5,0.5) to [out=-45,in=135] (3,0);
\draw[dashed, ultra thin] (3,0) to [out=225,in=45] (2.5,-0.5); 
\draw[dashed, ultra thin] (2.5,-0.5) to [out=135,in=-45] (2,0) ;

\draw[dashed, ultra thin] (3,0) to [out=45,in=135] (3.8,0);
\draw[dashed, ultra thin] (3.8,0) to [out=225,in=315] (3,0) ;
\draw[dashed, ultra thin] (3.8,0) to [out=45,in=135] (4.6,0) ;

\draw[ultra thin] (3.8,0) to [out=-45,in=-135] (4.6,0);

\draw[dashed, ultra thin] (0.3,0) circle [radius=0.1] ;
\draw[dashed, ultra thin] (2.5,0.6) circle [radius=0.1] ;
\draw[dashed, ultra thin] (2.5,-0.6) circle [radius=0.1] ;

\draw[fill] (.4,0) circle [radius=0.035] node[below=6pt]  {$(0,111)$};
\draw[fill] (1.2,0) circle [radius=0.035] node[below=6pt]  {$(0,011)$};
\draw[fill] (2,0) circle [radius=0.035] node[below=6pt]  {$(0,001)$};
\draw[fill] (2.5,0.5) circle [radius=0.035] node[right=4pt]  {$(0,100)$};
\draw[fill] (2.5,-0.5) circle [radius=0.035] node[left=6pt]  {$(0,101)$};
\draw[fill] (3,0) circle [radius=0.035] node[below=6pt]  {$(0,000)$};
\draw[fill] (3.8,0) circle [radius=0.035] node[above=6pt]  {$(0,010)$};
\draw[fill] (4.6,0) circle [radius=0.035] node[above=6pt]  {$(0,110)$};
\draw[dashed, ultra thin] (0.4,-1.5)to [out=45,in=135] (1.2,-1.5) ;
\draw[dashed, ultra thin] (1.2,-1.5)to [out=225,in=315] (0.4,-1.5) ;
\draw[dashed, ultra thin] (1.2,-1.5)to [out=45,in=135] (2,-1.5)  ;
\draw[dashed, ultra thin] (2,-1.5)to [out=225,in=315] (1.2,-1.5) ;

\draw[dashed, ultra thin] (2,-1.5) to [out=45,in=225] (2.5,-1) ;         
\draw[dashed, ultra thin] (2.5,-1) to [out=-45,in=135] (3,-1.5);
\draw[dashed, ultra thin] (3,-1.5) to [out=225,in=45] (2.5,-2); 
\draw[dashed, ultra thin] (2.5,-2) to [out=135,in=-45] (2,-1.5) ;

\draw[dashed, ultra thin] (3,-1.5) to [out=45,in=135] (3.8,-1.5);
\draw[dashed, ultra thin] (3.8,-1.5) to [out=225,in=315] (3,-1.5) ;
\draw[dashed, ultra thin] (3.8,-1.5) to [out=45,in=135] (4.6,-1.5) ;

\draw[ultra thin] (3.8,-1.5) to [out=-45,in=-135] (4.6,-1.5);

\draw[dashed, ultra thin] (0.3,-1.5) circle [radius=0.1] ;
\draw[dashed, ultra thin] (2.5,-0.9) circle [radius=0.1] ;
\draw[dashed, ultra thin] (2.5,-2.1) circle [radius=0.1] ;

\draw[fill] (.4,-1.5) circle [radius=0.035] node[below=6pt]  {$(1,111)$};
\draw[fill] (1.2,-1.5) circle [radius=0.035] node[below=6pt]  {$(1,011)$};
\draw[fill] (2,-1.5) circle [radius=0.035] node[below=6pt]  {$(1,001)$};
\draw[fill] (2.5,-1) circle [radius=0.035] node[right=4pt]  {$(1,100)$};
\draw[fill] (2.5,-2) circle [radius=0.035] node[left=6pt]  {$(1,101)$};
\draw[fill] (3,-1.5) circle [radius=0.035] node[below=6pt]  {$(1,000)$};
\draw[fill] (3.8,-1.5) circle [radius=0.035] node[below=6pt]  {$(1,010)$};
\draw[fill] (4.6,-1.5) circle [radius=0.035] node[below=6pt]  {$(1,110)$};
\draw[ultra thin] (4.6,0) to [out=-45,in=45] (4.6,-1.5);
\draw[ultra thin] (4.6,-1.5) to [out=135,in=-135] (4.6,0);
\draw[] (2.5,-2.5) node[]{$\Gamma_1\textcircled{g}\Gamma_3$};
\end{tikzpicture}
\caption{We choose $u_0 = 110.$ So the graph $\Gamma_1\textcircled{g}\Gamma_3$ is still normal $2$ sheeted cyclic covering of the graph $\Gamma_3$ but it is not isomorphic to $\Gamma_4$. The $2$ sheets of $\Gamma_3$ are copies of the spanning subgraph $\Gamma_3'$ (dashed lines).  Black  continuous lines/loops in $\Gamma_1\textcircled{g}\Gamma_3$ are lifts of $e_a$ or $e_b$.}
\label{fig2}
\end{figure}\\
 {\rem In the Figure 2,3 and 4 one can observe that the number of edges in the graph $\Gamma_n$ and the number of edges which are lifts of $e_a$ and $e_b$ in $\Gamma_n \textcircled{g} \Gamma_r$ are the same. If we imagine that the whole sheet as a single vertex along with the lift edges interestingly we get the graph $\Gamma_n$. Hence we say the following:
 \begin{enumerate} 
 \item If $r$ is even then the edges which are lifts of $e_a$ and $e_b$ in $\Gamma_n \textcircled{g} \Gamma_r$ are $a$ and $b$ edges in $\Gamma_n$ respectively.
 \item If $r$ is odd then the edges which are lifts of $e_a$ and $e_b$ in $\Gamma_n \textcircled{g} \Gamma_r$ are $b$ and $a$ edges in $\Gamma_n$ respectively.
 \item The edges which are not lifts of $e_a$ and $e_b$ in $\Gamma_n \textcircled{g} \Gamma_r$ are  edges in $\Gamma_r,$ so we call them sheet edges.
 \end{enumerate}}
Note that above definition of $\Gamma_n \textcircled{g} \Gamma_r$ depends on the spanning subgraph of the second graph. If a different choice of spanning subgraph is taken, then the graph $\Gamma_n \textcircled{g} \Gamma_r$ need not be Schreier graph of the Basilica group.\\
{\exa  If we choose a vertex of $\Gamma_r$ other than $0^r$ as $u_0$ and then corresponding spanning subgraph 
$\Gamma_r'$  with edge set, $E(\Gamma_r)\backslash \{e_a,e_b\}$ where $e_a$ and $e_b$ are the edges of the vertex $u_0$ with labels $a$ and $b$ respectively near $u_0,$ then the graph $\Gamma_n \textcircled{g} \Gamma_r$ need not be Schreier graph of the Basilica group. See Figure 2.}\\[7pt]
The following result will show that the graph $\Gamma_n \textcircled{g} \Gamma_r$ is actually a Schreier graph under the action of Basilica group on the set $X^{n+r}.$\\
{\Pro Let $n , r \geq 1$, then the following holds:
\begin{enumerate}
\item  The graph $\Gamma_{n+r} \simeq \Gamma_n\textcircled{g} \Gamma_r \simeq \Gamma_r \textcircled{g}   \Gamma_n $ is a Schreier graph under the action of Basilica group on the set $X^{n+r}.$
\item  $\Gamma_{n+r}$ is an unramified, $2^n$ sheeted graph covering of $\Gamma_r.$
\item  $\Gamma_n \textcircled{g} \Gamma_r$ is non normal covering of the graph $\Gamma_r$ if $n > 1.$
\item $\Gamma_{1+r}$ is $2$ sheeted normal covering of the graph $\Gamma_r.$
\end{enumerate} 
\begin{proof}
Let $n,r \geq 1.$
\begin{enumerate}
\item  Let us define the map $f: \Gamma_n \textcircled{g} \Gamma_r \rightarrow \Gamma_{r+n} $ by $f(v,u) = uv,$ where $v \in X^n$ and $u \in X^r.$ To show  $f$ is adjacency preserving map, notice that there are two types of edges in $\Gamma_n\textcircled{g}\Gamma_r,$ first type which contains those edges which are lifts of $e_a$ or $e_b$, i.e. edges which are described in the equation (3) or (4) and second type contains edges which are not lifts of $e_a$ or $e_b$, i.e. edges which are described in the equation (5).\\
{\bf Edges of the first type:}\\
Note that lift of $e_a$ connect $(v,u_0)$ to $(v',u_0^a)$ and lift of $e_b$ connect  $(v,u_0)$ to $(v'',u_0^b)$ where
$$ v' = \left\{              
                    \begin{array}{ll}           
                   v^a & \mathrm{if ~r ~is~ even} \\
                   v^b & \mathrm{if ~r~is~odd} \\                 
                    \end{array}       
                    \right. $$
$$ v'' = \left\{              
                    \begin{array}{ll}           
                   v^b & \mathrm{if~r~is~even} \\
                   v^a &  \mathrm{if~r~is~odd} \\                 
                    \end{array}       
                    \right. $$
Let given edge be a lift of $e_a.$ In other words, let $(v,u_0)$ be adjacent to $(v',u_0^a)$ in $\Gamma_n\textcircled{g}\Gamma_r,$ i.e.
$$ \Leftrightarrow {\bf Rot}_{\Gamma_n\textcircled{g}\Gamma_r}((v,u_0),a) = ((v',u_0^a),\alpha^{-1}) $$  
Then by action of the Basilica group $B$ on $X^{n+r},$
$$ \Leftrightarrow a(u_0v) = \left\{              
                    \begin{array}{ll}           
                   u_0^av^a & \mathrm{if~r~is~even} \\
                   u_0^av^b &\mathrm{if~r~is~odd} \\                 
                    \end{array}       
                    \right. $$
$$ \Leftrightarrow a(u_0v) = u_0^av',  v' = \left\{              
                    \begin{array}{ll}           
                   v^a & \mathrm{if~r~is~even} \\
                   v^b & \mathrm{if~r~is~odd} \\              
                    \end{array}       
                    \right.$$    
$$ \Leftrightarrow {\bf Rot}_{\Gamma_{n+r}}(u_0v,a) = (u_0^av'),  v' = \left\{              
                    \begin{array}{ll}           
                   v^a & \mathrm{if~r~is~even} \\
                   v^b & \mathrm{if~r~is~odd} \\                
                    \end{array}       
                    \right.$$
$\Leftrightarrow  u_0v$ is adjacent to $u_0^av'$ in $\Gamma_{n+r}$ which shows that $f$ is adjacency preserving map. The argument is similar if the given edge is a lift of $e_b.$\\
{\bf Edges of the second type:}\\
The edges of $\Gamma_n \textcircled{g} \Gamma_r$ which are not lifts of $e_a$ or $e_b$ connect $(v,u)$ to $(v,u^s)$ where $u, u^s \neq u_0.$ Suppose given edge is of this type.
$$ \Leftrightarrow {\bf Rot}_{\Gamma_n \textcircled{g} \Gamma_r}((v,u),s) = ((v,u^s),s^{-1}), s \in \{a^{\pm 1},b^{\pm 1}\}, u, u^s \neq u_0.$$
As $u \neq u_0 = 0^r$ means $u$ has at least one $``1"$ as an alphabet, therefore by definition of Basilica group elements we have 
$$ s(uv) = s(u)v = u^sv, s \in \{a^{\pm 1},b^{\pm 1}\} $$
Hence $uv$ is adjacent to $u^sv$ in $\Gamma_{n+r}$ which gives $f$ is adjacency preserving map. By definition of $f$ it is not difficult to that $f$ is bijection map.
Therefore $f^{-1}$ exists and by reverse implication it shows that $f^{-1}$ is also adjacency preserving map.
As $f$ bijection so as $f^{-1}.$ Thus $f$ and $f^{-1}$ both are isomorphisms. 
\item By definition 2.1, $\Gamma_{n+r} \simeq \Gamma_n \textcircled{g} \Gamma_r$ contains $2^n$ sheets of the graph $\Gamma_r$. To show that $\Gamma_{n+r}$ is an unramified cover of the graph $\Gamma_r$, we define the map 
$\pi : \Gamma_{n+r} \rightarrow \Gamma_r$ by $\pi(uv) = u.$ We now show that $\pi$ is a covering map i.e. $\pi$ is sending neighborhoods of $\Gamma_{n+r}$ one-to-one and onto neighborhoods of $\Gamma_r$.\\
Suppose $u_0v$ is adjacent to $u_0^sv'$ in $\Gamma_{n+r}$ where $s \in \{ a,b\}$.\\
By definition of $\pi$, $\pi(u_0v) = u_0$ and $\pi(u_0^sv') = u_0^s$ but $u_0^s = s(u_0)$
which means $u_0$ is adjacent to $s(u_0) = u_0^s$ in $\Gamma_r$ so $\pi(u_0v)$ is adjacent to $\pi(u_0^sv')$ in $\Gamma_r.$\\
Let $${\bf N}_{\Gamma_{n+r}}(u_0v) = \{u_0^av', u_0^bv'', u_0^{a^{-1}}v, u_0^{b^{-1}}v \},$$
then 
$$ \pi({\bf N}_{\Gamma_{n+r}}(u_0v)) = \pi(\{u_0^av', u_0^bv'', u_0^{a^{-1}}v, u_0^{b^{-1}}v \} )$$
$$ = \{u_0^a, u_0^b, u_0^{a^{-1}}, u_0^{b^{-1}} \} = {\bf N}_{\Gamma_{r}}(u_0)$$
Now suppose $uv$ is adjacent to $u^sv$ in $\Gamma_{n+r}$ where $s \in \{ a^{\pm 1}, b^{\pm 1}\}$ and $u, u^s \neq u_0.$\\
By definition of $\pi$, $\pi(uv) = u$ and $\pi(u^sv) = u^s$ but $u^s = s(u)$ in $\Gamma_r$ gives $u$ is adjacent to $u^s$ in $\Gamma_r.$\\
Let $${\bf N}_{\Gamma_{n+r}}(uv) = \{u^av, u^bv, u^{a^{-1}}v, u^{b^{-1}}v \},$$
then 
$$ \pi({\bf N}_{\Gamma_{n+r}}(uv)) = \pi(\{u^av, u^bv, u^{a^{-1}}v, u^{b^{-1}}v \} )$$
$$ = \{u^a, u^b, u^{a^{-1}}, u^{b^{-1}} \} = {\bf N}_{\Gamma_{r}}(u)$$
$\implies \pi$ is a covering map. Thus $\Gamma_{n+r}$ is an unramified covering of $\Gamma_r.$
\item (Proof by contradiction) Assume that $\Gamma_{n+r}$ is normal covering of $\Gamma_r,$ where $n > 1$. The Galois group $\mathbb{G} = Gal(\Gamma_{n+r}|\Gamma_r) = Gal(\Gamma_n \textcircled{g} \Gamma_r|\Gamma_r)$ is generated by permutations corresponding to the Frobenius automorphism $\sigma$ at $e_a$ and $e_b$. See \cite{Ter}\\
i.e. $$\mathbb{G} = ~< \sigma(e_a), \sigma(e_b) >$$
The permutations $\sigma(e_a)$ and $\sigma(e_b)$ arises under the action of Basilica group elements $a$ and $b$ on the (ordered) set $X^n$ respectively when we consider $\Gamma_{n+r}$ as unramified covering of $\Gamma_r.$\\ As $\Gamma_{n+r}$ has $2^n$ sheets of $\Gamma_r,$ we indexed these sheets by the set $\{1,2,\cdots,2^n\}$\\$ \sim X^n.$ If two sheets say $i^{th}$ and $j^{th}$ are connected by the lift of say $e_a$ then we write $\sigma(e_a)(i) = j$ similarly we can write the permutation corresponds to $\sigma(e_b).$\\
 The permutations $\sigma(e_a), \sigma(e_b)$ can be written directly using action of Basilica group elements $a$ and $b$ on the set $X^n$ respectively i.e. suppose $v_i, v_j \in X^n,$ with $a(v_i) = v_j$ then we write $\widetilde{\sigma}(e_a)(i) = j.$  By Remark 2.1, $\sigma \equiv \widetilde{\sigma}.$ Therefore the product $\sigma(e_a)\cdot\sigma(e_b)$ comes from the action of $ba$ on $X^n.$ Using induction, it can be shown as in \cite{DAng} that the action of $ba$ on $X^n$ has order $2^n$. This shows that the Galois group $\mathbb{G}$ has an element ($\sigma(e_a)\cdot\sigma(e_b)$) which is different from the generators $\sigma(e_a)$ and $\sigma(e_b)$ and which has order $2^n.$ Hence Galois group $\mathbb{G}$ has order bigger that $2^n$  which is a contradiction.  
\item Suppose $n = 1,$ therefore $\Gamma_{1+r}$ is $2$ sheeted covering of $\Gamma_r$ so we first show that $|\mathbb{G}| = 2.$ Denote two sheets $\Gamma_r'$ of $\Gamma_{1+r}$ by $1$ and $2.$\\
Recall that the Galois group $ \mathbb{G} = Gal(\Gamma_{1+r}|\Gamma_r) =~ <\sigma(e_a), \sigma(e_b) >.$ See \cite{Ter}. By the definition of $\Gamma_{1+n} \simeq \Gamma_1 \textcircled{g} \Gamma_n$ we have the following:\\ If $r$ is even then $$ \sigma(e_a) = (1), \sigma(e_b) = (1~ 2)$$
and if $r$ is odd then $$ \sigma(e_a) = (1 ~2), \sigma(e_b) = (1)$$
Therefore $$ \mathbb{G} = ~< \sigma = (1~ 2) > , ~~\forall~~ r. $$
To show that $\Gamma_{1+r}$ is $2$ sheeted normal covering of $\Gamma_r$ it remains to show that $\pi \circ \sigma \equiv \pi$ where $\sigma \in  \mathbb{G}$ is the non identity Frobenius automorphism with $\sigma^2 = id.$ If $u$ is any vertex of $\Gamma_r$ so let $u0$ be any vertex of the first sheet and $u1$ be any vertex of the second sheet.
$$ (\pi \circ \sigma)(u0) = \pi(\sigma(u0))  = \pi(u1) = u $$
$$ (\pi \circ \sigma)(u1) = \pi(\sigma(u1))  = \pi(u0) = u,$$
and $$ \pi(u1) = u = \pi(u0), $$
Therefore $$ \pi \circ \sigma \equiv \pi.$$
\end{enumerate}
\end{proof}
{\cor The spectrum of  $\Gamma_r$ is contained in the spectrum of $\Gamma_{1+r}$. 
Moreover The zeta function $\zeta_{\Gamma_r}(t)^{-1}$ divides $\zeta_{\Gamma_{1+r}}(t)^{-1}.$ }
\begin{proof}
This is direct consequence of the proof of the Proposition 13.10, p.no. 110 given in \cite{Ter}. 
\end{proof}
\textnormal{ Let $\widetilde{Y}$ be the normal covering of $Y$. The Artin $L$ function associated to the representation $\rho$ of $\mathbb{G} = Gal(\widetilde{Y}|Y)$ can be defined by a product over prime cycles in $Y$ as 
$$ L(u,\rho,\widetilde{Y}|Y) = \prod_{[c]~ prime~ in~ Y} \det(I - \rho(Frob(\widetilde{C})) ~~u^{\nu(C)})^{-1}$$
where $\widetilde{C}$ denotes any lift of $C$ to $\widetilde{Y}$ and $Frob(\widetilde{C})$ denotes the Frobenius automorphism defined by $$Frob(\widetilde{C}) = ji^{-1},$$
if $\widetilde{C}$ starts on $\widetilde{Y}$-sheet labeled by $i \in \mathbb{G}$ and ends on $\widetilde{Y}$-sheet labeled by $j \in \mathbb{G}.$ As in Proposition 2.1 of \cite{Ded}, the adjacency matrix of $\widetilde{Y}$  can be block diagonalized where the blocks are of the form
\begin{equation}
A_{\rho} = \sum_{g \in \mathbb{G}} A(g) \otimes \rho(g),
\end{equation}
each taken $d_{\rho} = $ degree of $\rho$ times, $\rho \in \widehat{\mathbb{G}}$ irreducible representation of $\mathbb{G}$ and $A(g)$ is defined in following formula:\\ Suppose $Y$ has $m$ vertices. Define the $m\times m$ matrix $A(g)$ for $g \in \mathbb{G}$ by defining $i,j$ entry to be 
\begin{equation}
A(g)_{i,j} = \mathrm{~the~number~of~edges~in}~\widetilde{Y}\mathrm{~between}~(1,e)\mathrm{~to}~(j,g), 
\end{equation}
where $e$ denotes the identity in $\mathbb{G}.$\\
By setting $Q_{\rho} = Q \otimes I_{d_{\rho}},$ with $d_{\rho} = deg ~\rho,$ we have the following analogue of formula (2):
\begin{equation}
L(t,\rho,\widetilde{Y}|Y)^{-1} = (1-t^2)^{(r-1)d_{\rho}} \det(I - A_{\rho}t + Q_{\rho}t^2).
\end{equation} 
See \cite{Ter}.\\
Thus we have zeta functions of $\widetilde{Y}$ factors as follows
\begin{equation}
\zeta_{\widetilde{Y}}(t) = \prod_{\rho \in \widehat{\mathbb{G}}}  L(t, \rho, \widetilde{Y}|Y)^{d_{\rho}}
\end{equation}
See Corollary 18.11 of \cite{Ter}.\\
The matrices $A_1$and $A_{\sigma}$ are called as the Artinized adjacency matrices.} \\
\begin{figure}
\begin{tikzpicture}[scale=2]
\draw[ultra thin] (-1,1) node[above right] {$^b$} to [out=45,in=135] (0,1) node[above] {$^{b^{-1}}$} ;
\draw[ultra thin] (0,1) node[below left] {$_{b}$} to [out=225,in=315] (-1,1) node[below right] {$_{b^{-1}}$};
\draw[ultra thin] (0,1) node[above right] {$^a$} to [out=45,in=135] (1,1) node[above] {$^{a^{-1}}$} ;
\draw[ultra thin] (1,1) node[below left] {$_{a}$} to [out=225,in=315] (0,1) node[below right] {$_{a^{-1}}$};
\draw[ultra thin] (1,1) node[above right] {$^b$} to [out=45,in=135] (2,1) node[above] {$^{b^{-1}}$} ;
\draw[ultra thin] (2,1) node[below left] {$_{b}$} to [out=225,in=315] (1,1) node[below right] {$_{b^{-1}}$};

\draw[] (2.1,1)node[above right] {$^{a^{-1}}$} circle [radius=0.1] node[below right] {$_{a}$};
\draw[] (-1.1,1)node[above] {$^{a^{-1}}$} circle [radius=0.1] node[below left] {$_{a}$};
\draw[fill] (-1,1) circle [radius=0.035] node[below=4pt]  {${\bf 10}$};
\draw[fill] (0,1) circle [radius=0.035] node[below=4pt]  {${\bf 00}$};
\draw[fill] (1,1) circle [radius=0.035] node[below=4pt]  {${\bf 01}$};
\draw[fill] (2,1) circle [radius=0.035] node[below=4pt]  {${\bf 11}$};
\draw[] (0.5,0.5) node[]{$\Gamma_2$}; 
\draw[color=gray, dashed, ultra thick] (3,1) node[above right] {$^b$} to [out=45,in=135] (4,1) node[above] {$^{b^{-1}}$} ;
\draw[color=gray, dashed, ultra thick] (5,1) node[below left] {$_{a}$} to [out=225,in=315] (4,1) node[below right] {$_{a^{-1}}$};
\draw[color=gray, dashed, ultra thick] (5,1) node[above right] {$^b$} to [out=45,in=135] (6,1) node[above] {$^{b^{-1}}$} ;
\draw[color=gray, dashed, ultra thick] (6,1) node[below left] {$_{b}$} to [out=225,in=315] (5,1) node[below right] {$_{b^{-1}}$};

\draw[color=gray, dashed, ultra thick] (6.1,1)node[above right] {$^{a^{-1}}$} circle [radius=0.1] node[below right] {$_{a}$};
\draw[color=gray, dashed, ultra thick] (2.9,1)node[above] {$^{a^{-1}}$} circle [radius=0.1] node[below left] {$_{a}$};
\draw[fill] (3,1) circle [radius=0.035] node[below=4pt]  {${\bf 10}$};
\draw[fill] (4,1) circle [radius=0.035] node[below=4pt]  {${\bf 00}$};
\draw[fill] (5,1) circle [radius=0.035] node[below=4pt]  {${\bf 01}$};
\draw[fill] (6,1) circle [radius=0.035] node[below=4pt]  {${\bf 11}$};
\draw[] (4.5,0.5) node[]{$\Gamma_2'$}; 
\draw[ ultra thin] (2,0) node[above right] {$^b$} to [out=45,in=135] (3,0) node[above] {$^{b^{-1}}$} ;
\draw[  ultra thin] (3,0) node[below left=5pt] {$_{b}$} to [out=225,in=315] (2,0) node[below right=7pt] {$_{b^{-1}}$};
\draw[  ultra thin] (3.1,0)node[above right] {$^{a^{-1}}$} circle [radius=0.1] node[below right] {$_{a}$};
\draw[ ultra thin] (1.9,0)node[above] {$^{a^{-1}}$} circle [radius=0.1] node[below left] {$_{a}$};
\draw[fill] (2,0) circle [radius=0.035] node[below=2pt]  {${\bf 0}$};
\draw[fill] (3,0) circle [radius=0.035] node[below=2pt]  {${\bf 1}$};
\draw[] (2.5,-0.5) node[]{$\Gamma_1$}; 

\draw[color=gray, dashed, ultra thick] (0.4,-1.5) node[above right] {$^b$} to [out=45,in=135] (1.2,-1.5) node[above] {$^{b^{-1}}$} ;
\draw[color=gray, dashed, ultra thick] (1.2,-1.5) node[below left] {$_{b}$} to [out=225,in=315] (0.4,-1.5) node[below right] {$_{b^{-1}}$};
\draw[color=gray, dashed, ultra thick] (1.2,-1.5) node[above right] {$^a$} to [out=45,in=135] (2,-1.5) node[above] {$^{a^{-1}}$} ;
\draw[ultra thin] (2,-1.5) node[below left] {$_{a}$} to [out=225,in=315] (1.2,-1.5) node[below right] {$_{a^{-1}}$};

\draw[color=gray, dashed, ultra thick] (2,-1.5) node[above right=2pt] {$^{b^{-1}}$} to [out=45,in=225] (2.5,-1) node[below left] {$^b$};        
\draw[ ultra thin] (2.5,-1) node[below right] {$^{b^{-1}}$} to [out=-45,in=135] (3,-1.5) node[above left] {$^b$};
\draw[color=gray, dashed, ultra thick] (3,-1.5) node[below left] {$_{b^{-1}}$} to [out=225,in=45] (2.5,-2) node[above right] {$_{b}$}; 
\draw[ultra thin] (2.5,-2) node[above left] {$_{b^{-1}}$} to [out=135,in=-45] (2,-1.5) node[below right] {$_{b}$};

\draw[ultra thin] (3,-1.5) node[above right] {$^a$} to [out=45,in=135] (3.8,-1.5) node[above] {$^{a^{-1}}$} ;
\draw[color=gray, dashed, ultra thick] (3.8,-1.5) node[below left] {$_{a}$} to [out=225,in=315] (3,-1.5) node[below right] {$_{a^{-1}}$};
\draw[color=gray, dashed, ultra thick] (3.8,-1.5) node[above right] {$^b$} to [out=45,in=135] (4.6,-1.5) node[above] {$^{b^{-1}}$} ;
\draw[color=gray, dashed, ultra thick] (4.6,-1.5) node[below left] {$_{b}$} to [out=225,in=315] (3.8,-1.5) node[below right] {$_{b^{-1}}$};

\draw[color=gray, dashed, ultra thick] (4.7,-1.5) node[above right] {$^{a^{-1}}$} circle [radius=0.1] node[below right] {$_{a}$};
\draw[color=gray, dashed, ultra thick] (0.3,-1.5) node[above] {$^{a^{-1}}$} circle [radius=0.1] node[below left] {$_{a}$};
\draw[color=gray, dashed, ultra thick] (2.5,-0.9) node[above right=-1pt] {$^{a^{-1}}$} circle [radius=0.1] node[above left] {$_{a}$};
\draw[color=gray, dashed, ultra thick] (2.5,-2.1) node[below right=-1pt ] {$^{a^{-1}}$} circle [radius=0.1] node[below left] {$_{a}$};

\draw[fill] (.4,-1.5) circle [radius=0.035] node[below=6pt]  {${\bf 1}$};
\draw[fill] (1.2,-1.5) circle [radius=0.035] node[below=6pt]  {${\bf 2}$};
\draw[fill] (2,-1.5) circle [radius=0.035] node[below=6pt]  {${\bf 3}$};
\draw[fill] (2.5,-1) circle [radius=0.035] node[right=4pt]  {${\bf 4}$};
\draw[fill] (2.5,-2) circle [radius=0.035] node[left=6pt]  {${\bf 4'}$};
\draw[fill] (3,-1.5) circle [radius=0.035] node[below=6pt]  {${\bf 3'}$};
\draw[fill] (3.8,-1.5) circle [radius=0.035] node[below=6pt]  {${\bf 2'}$};
\draw[fill] (4.6,-1.5) circle [radius=0.035] node[below=6pt]  {${\bf 1'}$};
\draw[] (2.5,-2.5) node[]{$\Gamma_1 \textcircled{g} \Gamma_2$}; 
\draw[->] (2.5,-2.7) to (2.5,-3.6);
\draw[](2.6,-3.1)node[right] {$f$};
\draw[ultra thin] (0.4,-4.5) node[above right] {$^b$} to [out=45,in=135] (1.2,-4.5) node[above] {$^{b^{-1}}$} ;
\draw[ultra thin] (1.2,-4.5) node[below left] {$_{b}$} to [out=225,in=315] (0.4,-4.5) node[below right] {$_{b^{-1}}$};
\draw[ultra thin] (1.2,-4.5) node[above right] {$^a$} to [out=45,in=135] (2,-4.5) node[above] {$^{a^{-1}}$} ;
\draw[ultra thin] (2,-4.5) node[below left] {$_{a}$} to [out=225,in=315] (1.2,-4.5) node[below right] {$_{a^{-1}}$};

\draw[ultra thin] (2,-4.5) node[above right=2pt] {$^{b^{-1}}$} to [out=45,in=225] (2.5,-4) node[below left] {$^b$};         \draw[ultra thin] (2.5,-4) node[below right] {$^{b^{-1}}$} to [out=-45,in=135] (3,-4.5) node[above left] {$^b$};
\draw[ultra thin] (3,-4.5) node[below left] {$_{b^{-1}}$} to [out=225,in=45] (2.5,-5) node[above right] {$_{b}$}; 
\draw[ultra thin] (2.5,-5) node[above left] {$_{b^{-1}}$} to [out=135,in=-45] (2,-4.5) node[below right] {$_{b}$};

\draw[ultra thin] (3,-4.5) node[above right] {$^a$} to [out=45,in=135] (3.8,-4.5) node[above] {$^{a^{-1}}$} ;
\draw[ultra thin] (3.8,-4.5) node[below left] {$_{a}$} to [out=225,in=315] (3,-4.5) node[below right] {$_{a^{-1}}$};
\draw[ultra thin] (3.8,-4.5) node[above right] {$^b$} to [out=45,in=135] (4.6,-4.5) node[above] {$^{b^{-1}}$} ;
\draw[ultra thin] (4.6,-4.5) node[below left] {$_{b}$} to [out=225,in=315] (3.8,-4.5) node[below right] {$_{b^{-1}}$};

\draw[] (4.7,-4.5) node[above right] {$^{a^{-1}}$} circle [radius=0.1] node[below right] {$_{a}$};
\draw[] (0.3,-4.5) node[above] {$^{a^{-1}}$} circle [radius=0.1] node[below left] {$_{a}$};
\draw[] (2.5,-3.9) node[above right=-1pt] {$^{a^{-1}}$} circle [radius=0.1] node[above left] {$_{a}$};
\draw[] (2.5,-5.1) node[below right=-1pt ] {$^{a^{-1}}$} circle [radius=0.1] node[below left] {$_{a}$};

\draw[fill] (.4,-4.5) circle [radius=0.035] node[below=6pt]  {${\bf 110}$};
\draw[fill] (1.2,-4.5) circle [radius=0.035] node[below=6pt]  {${\bf 010}$};
\draw[fill] (2,-4.5) circle [radius=0.035] node[below=6pt]  {${\bf 000}$};
\draw[fill] (2.5,-4) circle [radius=0.035] node[right=4pt]  {${\bf 100}$};
\draw[fill] (2.5,-5) circle [radius=0.035] node[left=6pt]  {${\bf 101}$};
\draw[fill] (3,-4.5) circle [radius=0.035] node[below=6pt]  {${\bf 001}$};
\draw[fill] (3.8,-4.5) circle [radius=0.035] node[below=6pt]  {${\bf 011}$};
\draw[fill] (4.6,-4.5) circle [radius=0.035] node[below=6pt]  {${\bf 111}$};
\draw[] (2.5,-5.5) node[]{$\Gamma_3$}; 
\end{tikzpicture}
\caption{The graph $\Gamma_1\textcircled{g}\Gamma_2$ is a normal $2$ sheeted cyclic covering of the graph $\Gamma_2$. The $2$ sheets of $\Gamma_2$ are copies of the spanning subgraph $\Gamma_2'$ (dashed lines). The map $f$ is a graph isomorphism between $\Gamma_1\textcircled{g}\Gamma_2$ and $\Gamma_3$.}
\label{fig2}
\end{figure}
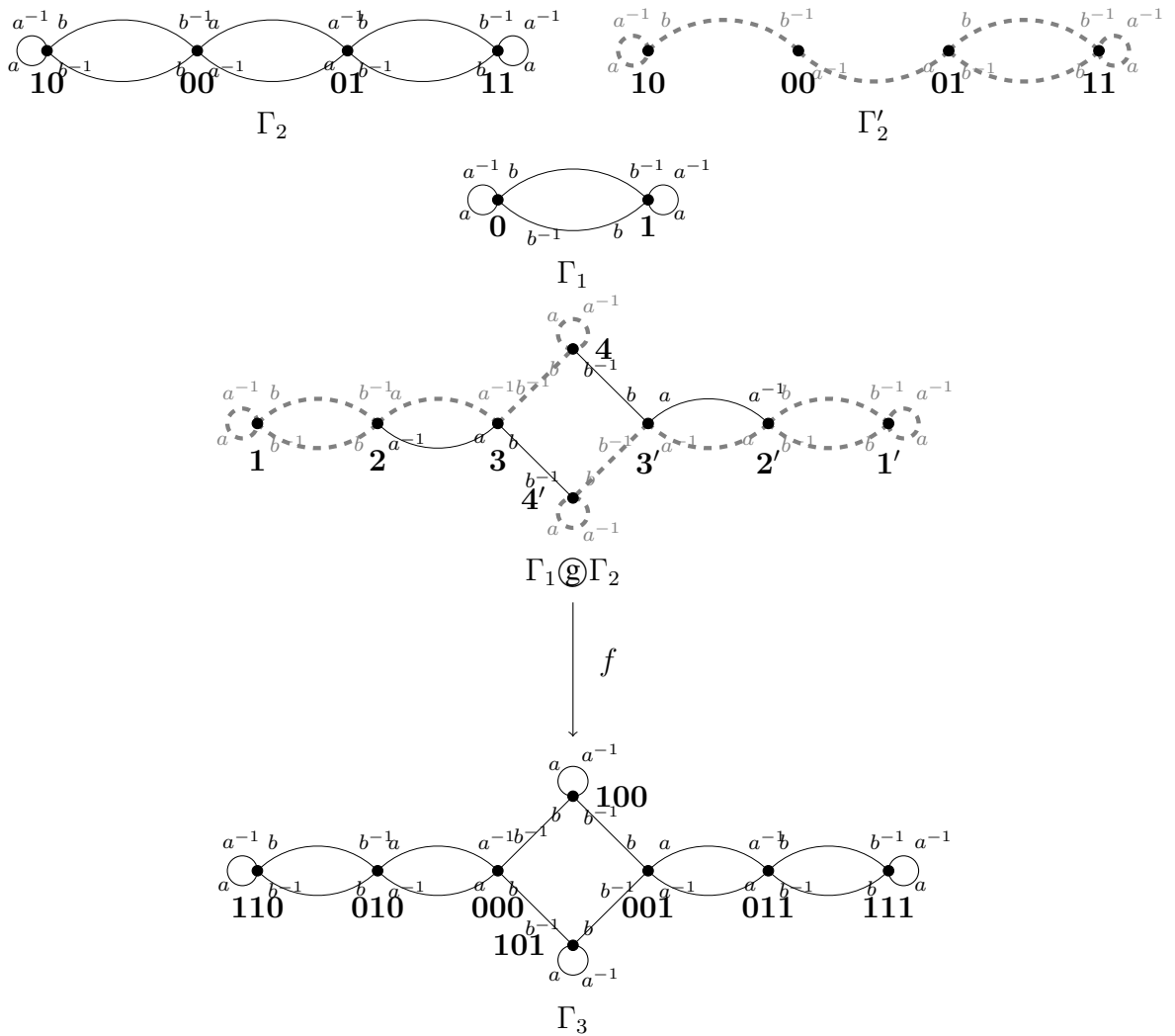\\
We label the vertices using the Table 1.
\begin{table}[h]
\begin{center}
\caption{Vertex labeling for $\Gamma_{2}$}
\begin{tabular}{|c|c|c|c|c|}
\hline
label & 1 & 2 & 3 & 4 \\
\hline
\hline
vertex & $11$  & $01$ & $00$ & $10$ \\
\hline
\end{tabular}
\end{center}
\end{table}\\
{\exa The covering $\Gamma_{3} \simeq \Gamma_{1} \textcircled{g}\Gamma_2 $ over $\Gamma_{2}$ given in the Figure 3 is a normal covering.}
We obtain spanning subgraph of the $Y = \Gamma_{2}$ by cutting edges $e_a = \{3,2\},e_b = \{3,4\}.$ This gives dashed line spanning subgraph $\Gamma_2'$ of the $\Gamma_{2}$. So we draw the covering graph $\widetilde{Y} = \Gamma_3 \simeq \Gamma_1\textcircled{g}\Gamma_2$ by placing $2$ sheets of spanning subgraph $\Gamma_2'$ of $Y = \Gamma_2$ and labeling each sheet given in Table 2.
Connections between $2$ sheets in the cover graph $\widetilde{Y} = \Gamma_3$ are given in the Table 3. In this case the representations of cyclic Galois group $\mathbb{G} = \{ 1, \sigma \}$ are the trivial representation $\rho_0 = 1$ and the representation $\rho$ defined by $\rho(1) = 1, \rho(\sigma) = -1.$ So $Q_{\rho} = 2 I_4$. There are two cases. \\
\begin{table}[h]
\caption{Notations for sheets}
\begin{center}
\begin{tabular}{|c|c|c|}
\hline
Vertex set of $\Gamma_1\textcircled{g}\Gamma_2$ & Vertex set of $\Gamma_3$ & Group element\\
\hline
\hline
$\{1,2,3,4\} \sim \{(0,11), (0,01),(0,00),(0,10)\}$ & $\{110, 010,000,100 \}$ & $1 = id$\\
$\{1',2',3',4'\} \sim \{(1,11), (1,01),(1,00),(1,10)\}$ & $\{111, 011,001,101 \}$ & $\sigma \neq id$\\
\hline
\end{tabular}
\end{center}
  \caption{Connections between sheet 1 and 2 in $\Gamma_3 \simeq \Gamma_1 \textcircled{g} \Gamma_2$}
 \begin{center}
 \begin{tabular}{|c|c|}
 \hline
  Vertex & adjacent vertices in $\Gamma_3 \simeq \Gamma_1 \textcircled{g} \Gamma_2$\\
  \hline
  \hline
  $1$ & $1,2$\\
  $2$ & $1,3$\\
  $3$ & $2,4,4'$\\
  $4$ & $3,4,3'$\\
 \hline
 \end{tabular}
 \end{center}
 \end{table}\\
${\bf Case ~~1}$ {\bf The trivial representation} $\rho_0 = {\bf 1}$\\
Here $A_1 = A(1) + A(\sigma) = A_{Y}$ where 
$$ A(1) = \begin{pmatrix}
 2 & 2 & 0 & 0 \\
2 & 0 & 2 & 0 \\
0 & 2 & 0 & 1 \\
0 & 0 & 1 & 2 \\
\end{pmatrix}, A(\sigma) = \begin{pmatrix}
0 & 0 & 0 & 0\\
0 & 0 & 0 & 0\\
0 & 0 & 0 & 1\\
0 & 0 & 1 & 0\\
\end{pmatrix} $$
and  $A_1$ is the adjacency matrix  $A_{Y}.$\\
${\bf Case ~~2}$ {\bf The representation} $\rho$\\
$$ A_{\sigma} = A(1) - A(\sigma) =  \begin{pmatrix}
 2 & 2 & 0 & 0 \\
2 & 0 & 2 & 0 \\
0 & 2 & 0 & 0 \\
0 & 0 & 0 & 2 \\
\end{pmatrix}$$
{\large Reciprocals of $L$ functions for $\Gamma_{3}|\Gamma_{2} $.}\\
1) For $A_1$
$$ \zeta_{\Gamma_2}(t)^{-1} = L(t,A_1,\Gamma_{3}|\Gamma_{2} )^{-1} = (1-t^2)^4 (t-1) (3 t-1) \left(3 t^2+1\right) \left(9 t^4-2 t^2+1\right).$$
2) For $A_{\sigma}$
$$ L(t,A_{\sigma},\widetilde{Y}|Y)^{-1} = (1-t^2)^4\left(3 t^2-2 t+1\right) \left(27 t^6-18 t^5+3 t^4-4 t^3+t^2-2 t+1\right) $$
 By equation (9) we have
$$\zeta_{\Gamma_{3}}(t)^{-1} = L(t,A_1,\widetilde{Y}|Y)^{-1}L(t,A_{\sigma},\widetilde{Y}|Y)^{-1}$$ 
\begin{eqnarray}
= (1-t^2)^8 (t-1) (3 t-1) \left(3 t^2+1\right) \left(3 t^2-2 t+1\right) \left(9 t^4-2 t^2+1\right) \nonumber\\ \left(27 t^6-18 t^5+3 t^4-4 t^3+t^2-2 t+1\right) \nonumber
\end{eqnarray}
{\exa The coverings $\Gamma_{5} \simeq \Gamma_{3} \textcircled{g}\Gamma_2 $ over $\Gamma_{2}$ given in the Figure 4 is non normal covering.}
We follow the labeling the vertices given in Table 1. We obtain spanning subgraph of the $\Gamma_{2}$ by cutting edges $e_a = \{3,2\},e_b = \{3,4\}.$ 
This gives dashed line spanning subgraph $\Gamma_2'$ of the $\Gamma_{2}$. So we draw the covering graph $ \Gamma_5 \simeq \Gamma_3\textcircled{g}\Gamma_2$ by placing $8$ sheets of spanning subgraph $\Gamma_2'$ of $\Gamma_2$
and labeling each with  notations given in Table 4 and connections between first sheet with other sheets are given in Table 5. 
\begin{figure}
\begin{tikzpicture}[scale=1.5]
\draw[ultra thin] (-1,1) node[above right] {$^b$} to [out=45,in=135] (0,1) node[above] {$^{b^{-1}}$} ;
\draw[ultra thin] (0,1) node[below left] {$_{b}$} to [out=225,in=315] (-1,1) node[below right] {$_{b^{-1}}$};
\draw[ultra thin] (0,1) node[above right] {$^a$} to [out=45,in=135] (1,1) node[above] {$^{a^{-1}}$} ;
\draw[ultra thin] (1,1) node[below left] {$_{a}$} to [out=225,in=315] (0,1) node[below right] {$_{a^{-1}}$};
\draw[ultra thin] (1,1) node[above right] {$^b$} to [out=45,in=135] (2,1) node[above] {$^{b^{-1}}$} ;
\draw[ultra thin] (2,1) node[below left] {$_{b}$} to [out=225,in=315] (1,1) node[below right] {$_{b^{-1}}$};

\draw[] (2.1,1)node[above right] {$^{a^{-1}}$} circle [radius=0.1] node[below right] {$_{a}$};
\draw[] (-1.1,1)node[above] {$^{a^{-1}}$} circle [radius=0.1] node[below left] {$_{a}$};
\draw[fill] (-1,1) circle [radius=0.035] node[below=4pt]  {${\bf 10}$};
\draw[fill] (0,1) circle [radius=0.035] node[below=4pt]  {${\bf 00}$};
\draw[fill] (1,1) circle [radius=0.035] node[below=4pt]  {${\bf 01}$};
\draw[fill] (2,1) circle [radius=0.035] node[below=4pt]  {${\bf 11}$};
\draw[] (0.5,0.5) node[]{$\Gamma_2$}; 
\hspace{.5cm}
\draw[color=gray, dashed, ultra thick] (3,1) node[above right] {$^b$} to [out=45,in=135] (4,1) node[above] {$^{b^{-1}}$} ;
\draw[color=gray, dashed, ultra thick] (5,1) node[below left] {$_{a}$} to [out=225,in=315] (4,1) node[below right] {$_{a^{-1}}$};
\draw[color=gray, dashed, ultra thick] (5,1) node[above right] {$^b$} to [out=45,in=135] (6,1) node[above] {$^{b^{-1}}$} ;
\draw[color=gray, dashed, ultra thick] (6,1) node[below left] {$_{b}$} to [out=225,in=315] (5,1) node[below right] {$_{b^{-1}}$};

\draw[color=gray, dashed, ultra thick] (6.1,1)node[above right] {$^{a^{-1}}$} circle [radius=0.1] node[below right] {$_{a}$};
\draw[color=gray, dashed, ultra thick] (2.9,1)node[above] {$^{a^{-1}}$} circle [radius=0.1] node[below left] {$_{a}$};
\draw[fill] (3,1) circle [radius=0.035] node[below=4pt]  {${\bf 10}$};
\draw[fill] (4,1) circle [radius=0.035] node[below=4pt]  {${\bf 00}$};
\draw[fill] (5,1) circle [radius=0.035] node[below=4pt]  {${\bf 01}$};
\draw[fill] (6,1) circle [radius=0.035] node[below=4pt]  {${\bf 11}$};
\draw[] (4.5,0.5) node[]{$\Gamma_2'$}; 
\end{tikzpicture}\\
\vspace{.5cm}

\begin{center}
\begin{tikzpicture}[scale=2]

\draw[ultra thin] (0.4,-4.5) node[above right] {$^b$} to [out=45,in=135] (1.2,-4.5) node[above] {$^{b^{-1}}$} ;
\draw[ultra thin] (1.2,-4.5) node[below left] {$_{b}$} to [out=225,in=315] (0.4,-4.5) node[below right] {$_{b^{-1}}$};
\draw[ultra thin] (1.2,-4.5) node[above right] {$^a$} to [out=45,in=135] (2,-4.5) node[above] {$^{a^{-1}}$} ;
\draw[ultra thin] (2,-4.5) node[below left] {$_{a}$} to [out=225,in=315] (1.2,-4.5) node[below right] {$_{a^{-1}}$};

\draw[ultra thin] (2,-4.5) node[above right=2pt] {$^{b^{-1}}$} to [out=45,in=225] (2.5,-4) node[below left] {$^b$};         \draw[ultra thin] (2.5,-4) node[below right] {$^{b^{-1}}$} to [out=-45,in=135] (3,-4.5) node[above left] {$^b$};
\draw[ultra thin] (3,-4.5) node[below left] {$_{b^{-1}}$} to [out=225,in=45] (2.5,-5) node[above right] {$_{b}$}; 
\draw[ultra thin] (2.5,-5) node[above left] {$_{b^{-1}}$} to [out=135,in=-45] (2,-4.5) node[below right] {$_{b}$};

\draw[ultra thin] (3,-4.5) node[above right] {$^a$} to [out=45,in=135] (3.8,-4.5) node[above] {$^{a^{-1}}$} ;
\draw[ultra thin] (3.8,-4.5) node[below left] {$_{a}$} to [out=225,in=315] (3,-4.5) node[below right] {$_{a^{-1}}$};
\draw[ultra thin] (3.8,-4.5) node[above right] {$^b$} to [out=45,in=135] (4.6,-4.5) node[above] {$^{b^{-1}}$} ;
\draw[ultra thin] (4.6,-4.5) node[below left] {$_{b}$} to [out=225,in=315] (3.8,-4.5) node[below right] {$_{b^{-1}}$};

\draw[] (4.7,-4.5) node[above right] {$^{a^{-1}}$} circle [radius=0.1] node[below right] {$_{a}$};
\draw[] (0.3,-4.5) node[above] {$^{a^{-1}}$} circle [radius=0.1] node[below left] {$_{a}$};
\draw[] (2.5,-3.9) node[above right=-1pt] {$^{a^{-1}}$} circle [radius=0.1] node[above left] {$_{a}$};
\draw[] (2.5,-5.1) node[below right=-1pt ] {$^{a^{-1}}$} circle [radius=0.1] node[below left] {$_{a}$};

\draw[fill] (.4,-4.5) circle [radius=0.035] node[below=6pt]  {${\bf 110}$};
\draw[fill] (1.2,-4.5) circle [radius=0.035] node[below=6pt]  {${\bf 010}$};
\draw[fill] (2,-4.5) circle [radius=0.035] node[below=6pt]  {${\bf 000}$};
\draw[fill] (2.5,-4) circle [radius=0.035] node[right=4pt]  {${\bf 100}$};
\draw[fill] (2.5,-5) circle [radius=0.035] node[left=6pt]  {${\bf 101}$};
\draw[fill] (3,-4.5) circle [radius=0.035] node[below=6pt]  {${\bf 001}$};
\draw[fill] (3.8,-4.5) circle [radius=0.035] node[below=6pt]  {${\bf 011}$};
\draw[fill] (4.6,-4.5) circle [radius=0.035] node[below=6pt]  {${\bf 111}$};
\draw[] (2.5,-5.5) node[]{$\Gamma_3$}; 
\end{tikzpicture}\\
\vspace{.5cm}
\begin{tikzpicture}[scale=1.5]


\draw[fill] (0,1) circle [radius=0.025] node[below=6pt]  {$1$};
\draw[fill] (1,1) circle [radius=0.025] node[below=6pt]  {$2$};
\draw[fill] (2,1) circle [radius=0.025]node[below=6pt]  {$3$};
\draw[fill] (2.5,1.5) circle [radius=0.025]node[left=5pt]  {$4$};
\draw[fill] (2.5,0.5) circle [radius=0.025]node[left=5pt]  {$4'$};
\draw[fill] (3,1) circle [radius=0.025]node[below=6pt]  {$3'$};
\draw[fill] (3.5,1.5) circle [radius=0.025]node[left=5pt]  {$2'$};
\draw[fill] (3.5,2.3) circle [radius=0.025]node[left=5pt]  {$1'$};
\draw[fill] (3.5,0.5) circle [radius=0.025]node[left=5pt]  {$2''$};
\draw[fill] (3.5,-0.3) circle [radius=0.025]node[left=5pt]  {$1''$};
\draw[fill] (4,1) circle [radius=0.025]node[left=5pt]  {$3''$};
\draw[fill] (4.3,1.6) circle [radius=0.025]node[left=5pt]  {$4''$};
\draw[fill] (5,1.8) circle [radius=0.025]node[below=5pt]  {$3'''$};
\draw[fill] (5,2.6) circle [radius=0.025]node[left=5pt]  {$2'''$};
\draw[fill] (5,3.4) circle [radius=0.025]node[left=5pt]  {$1'''$};
\draw[fill] (5.7,1.6) circle [radius=0.025]node[right=5pt]  {$4'''$};
\draw[fill] (5,0.2) circle [radius=0.025]node[right=5pt]  {$3^{(5)}$};
\draw[fill] (4.3,0.4) circle [radius=0.025]node[right=5pt]  {$4^{(5)}$};
\draw[fill] (5,-0.6) circle [radius=0.025]node[right=5pt]  {$2^{(5)}$};
\draw[fill] (5,-1.4) circle [radius=0.025]node[right=5pt]  {$1^{(5)}$};
\draw[fill] (5.7,0.4) circle [radius=0.025]node[right=5pt]  {$4^{(6)}$};
\draw[fill] (6,1) circle [radius=0.025]node[right=5pt]  {$3^{(6)}$};
\draw[fill] (6.5,1.5) circle [radius=0.025]node[right=5pt]  {$2^{(6)}$};
\draw[fill] (6.5,2.3) circle [radius=0.025]node[right=5pt]  {$1^{(6)}$};
\draw[fill] (6.5,0.5) circle [radius=0.025]node[right=5pt]  {$2^{(7)}$};
\draw[fill] (6.5,-0.3) circle [radius=0.025]node[right=5pt]  {$1^{(7)}$};
\draw[fill] (7,1) circle [radius=0.025]node[right=5pt]  {$3^{(7)}$};
\draw[fill] (7.5,1.5) circle [radius=0.025]node[right=5pt]  {$4^{(7)}$};
\draw[fill] (7.5,0.5) circle [radius=0.025]node[right=5pt]  {$4^{(8)}$};
\draw[fill] (8,1) circle [radius=0.025]node[above=3pt]  {$3^{(8)}$};
\draw[fill] (9,1) circle [radius=0.025]node[above=3pt]  {$2^{(8)}$};
\draw[fill] (10,1) circle [radius=0.025]node[above=3pt]  {$1^{(8)}$};

\draw[dashed, ultra thin] (0,1)  to [out=45,in=135] (1,1) to [out=45,in=135](2,1) to [out = 45 , in =225 ]  (2.5,1.5);

\draw[dashed, ultra thin] (0,1)  to [out=-45,in=225] (1,1)  ;

\draw[dashed, ultra thin](2.5,0.5) to [out = 45, in =225 ](3,1)to [out = 45 , in = 225 ](3.5,1.5) to [out = 45 , in = -45](3.5,2.3) to [out =225 , in = 135](3.5,1.5); 

\draw[dashed, ultra thin](4.3,1.6) to (4,1) to (3.5,0.5) to [out =225 , in = 135](3.5,-0.3) to  [out =45 , in =-45](3.5,0.5) ; 

\draw[dashed, ultra thin](5.7,1.6) to (5,1.8);
\draw[dashed, ultra thin](5,2.6) to  [out =135 , in =225](5,3.4) to [out =-45 , in = 45](5,2.6) to [out =-45 , in = 45](5,1.8) ; 

\draw[dashed, ultra thin](4.3,0.4) to (5,0.2);
\draw[dashed, ultra thin](5,-0.6) to  [out =225 , in =135](5,-1.4) to [out =45 , in = -45](5,-0.6) to [out =45 , in = -45](5,0.2) ; 

\draw[dashed, ultra thin](5.7,0.4) to (6,1) to (6.5,1.5) to [out =135 , in = 225](6.5,2.3) to  [out =-45 , in =45](6.5,1.5) ; 

\draw[dashed, ultra thin](7.5,1.5) to [out = 225, in =45 ](7,1)to [out =225 , in =45 ](6.5,0.5) to [out =225 , in =135](6.5,-0.3) to [out =45 , in =-45](6.5,0.5); 

\draw[dashed, ultra thin]  (9,1) to [out=45,in=135](10,1);

\draw[dashed, ultra thin] (10,1)  to [out=225,in=-45] (9,1) to [out=225,in=-45](8,1) to [out = 225 , in =-135 ]  (7.5,0.5);


\draw[dashed, ultra thin] (-.16,1) circle [radius=0.16];
\draw[dashed, ultra thin] (2.5,1.66) circle [radius=0.16];
\draw[dashed, ultra thin] (2.5,0.34) circle [radius=0.16];
\draw[dashed, ultra thin] (3.5,2.46) circle [radius=0.16];
\draw[dashed, ultra thin] (3.5,-0.46) circle [radius=0.16];
\draw[dashed, ultra thin] (4.3,1.76) circle [radius=0.16];
\draw[dashed, ultra thin] (4.3,0.24) circle [radius=0.16];
\draw[dashed, ultra thin] (5.7,1.76) circle [radius=0.16];
\draw[dashed, ultra thin] (5.7,0.24) circle [radius=0.16];
\draw[dashed, ultra thin] (5,3.56) circle [radius=0.16];
\draw[dashed, ultra thin] (5,-1.56) circle [radius=0.16];
\draw[dashed, ultra thin] (6.5,2.46) circle [radius=0.16];
\draw[dashed, ultra thin] (6.5,-0.46) circle [radius=0.16];
\draw[dashed, ultra thin] (7.5,1.66) circle [radius=0.16];
\draw[dashed, ultra thin] (7.5,0.36) circle [radius=0.16];
\draw[dashed, ultra thin] (10.16,1) circle [radius=0.16];
\draw[] (1,1) to [out=-45 , in= 225](2,1) to (2.5,0.5);
\draw[] (2.5,1.5) to (3,1) to (3.5,0.5);
\draw[] (3.5,1.5) to (4,1) to (4.3,0.4);
\draw[] (4.3,1.6) to (5,1.8) to [out=135 ,in=225 ](5,2.6);
\draw[] (5.7,0.4) to (5,0.2) to [out=225 ,in=135 ](5,-0.6);
\draw[] (6.5,0.5) to (6,1) to (5.7,1.6);
\draw[] (7.5,0.5) to (7,1) to (6.5,1.5);
\draw[] (9,1) to [out=135 , in= 45](8,1) to (7.5,1.5);
\draw[] (5,-2.2) node[]{$\Gamma_3 \textcircled{g}\Gamma_2 \simeq \Gamma_5$}; 
\end{tikzpicture}
\end{center}
\caption{The graph $\Gamma_3\textcircled{g}\Gamma_2$ is a non normal $8$ sheeted  covering of the graph $\Gamma_2$. The $8$ sheets of $\Gamma_2$ are copies of the spanning subgraph $\Gamma_2'$ (dashed lines). The graph $\Gamma_2\textcircled{g}\Gamma_3$ is isomorphic to the graph $\Gamma_5$.}
\label{}
\end{figure}
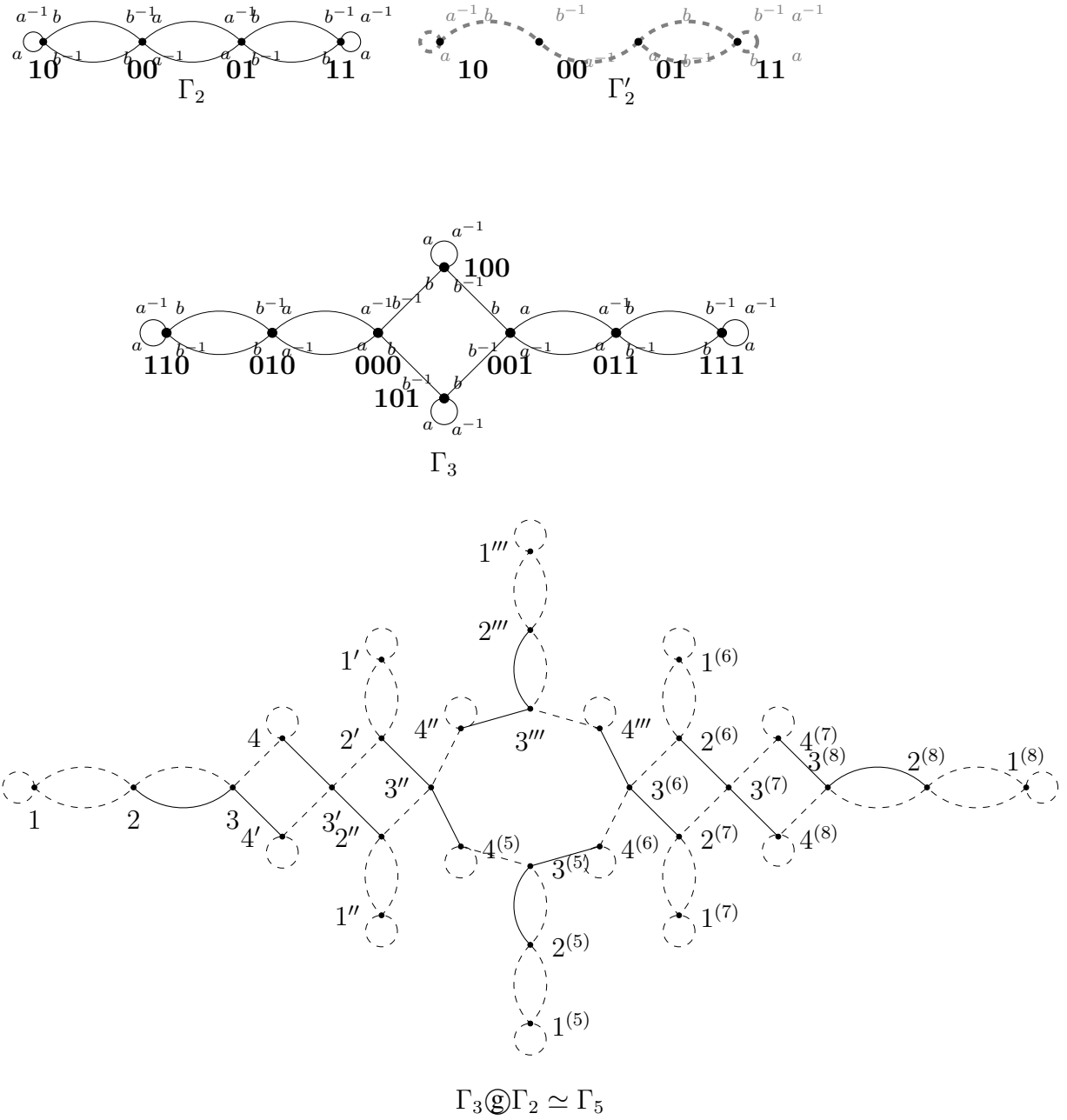
\begin{table}[h]
\caption{Notations for sheets}
\begin{center}
\begin{tabular}{|c|c|c|}
\hline
Vertex set of $\Gamma_3\textcircled{g}\Gamma_2$ & Sheet index\\
\hline
\hline
$\{1,2,3,4\} \sim \{(110,11), (110,01),(110,00),(110,10)\}$ &  $1^{st}$\\
$\{1',2',3',4'\} \sim \{(010,11), (010,01),(010,00),(010,10)\}$ &  $2^{nd}$\\
$\{1'',2'',3'',4''\} \sim \{(000,11), (000,01),(000,00),(000,10)\}$ & $3^{rd}$\\
$\{1''',2''',3''',4'''\} \sim \{(100,11), (100,01),(100,00),(100,10)\}$ & $4^{th}$\\
$\{1^{(5)},2^{(5)},3^{(5)},4^{(5)}\} \sim \{(101,11), (101,01),(101,00),(101,10)\}$ & $5^{th}$\\
$\{1^{(6)},2^{(6)},3^{(6)},4^{(6)}\} \sim \{(001,11), (001,01),(001,00),(001,10)\}$ & $6^{th}$\\
$\{1^{(7)},2^{(7)},3^{(7)},4^{(7)}\} \sim \{(011,11), (011,01),(011,00),(011,10)\}$ & $7^{th}$\\
$\{1^{(8)},2^{(8)},3^{(8)},4^{(8)}\} \sim \{(111,11), (111,01),(111,00),(111,10)\}$ & $8^{th}$\\
\hline
\end{tabular}
\end{center}
\end{table}\\
\begin{table}[h]
\caption{Connections between sheet 1 and other sheets in $\Gamma_5 \simeq \Gamma_3 \textcircled{g} \Gamma_2$}
\begin{center}
\begin{tabular}{|c|c|}
\hline
Vertex & adjacent vertices in $\Gamma_5 \simeq \Gamma_3 \textcircled{g} \Gamma_2$\\
  \hline
  \hline
  $1$ & $1,2$\\
  $2$ & $1,3$\\
  $3$ & $2,4,4'$\\
  $4$ & $3,4,3'$\\
 \hline
 \end{tabular}
 \end{center}
 \end{table}\\
Here we assume that $\Gamma_3\textcircled{g}\Gamma_2 $ is normal covering of $\Gamma_2.$
$$\sigma(e_a) = (2 ~3)(6~ 7), \sigma(e_b) = (1~ 2)(3~ 5~ 6~ 4)(7~ 8)$$ and $$\sigma(e_a)\circ\sigma(e_b) = (1 ~3~ 5~ 7~ 8~ 6~ 4~ 2)$$
Hence $\mathbb{G} = Gal(\Gamma_3\textcircled{g}\Gamma_2|\Gamma_2) =  < \sigma(e_a), \sigma(e_b) > $ has an element of order $8.$ Thus $$|\mathbb{G}| > 8 $$ is a contradiction.
\begin{center}
$\implies \Gamma_3\textcircled{g}\Gamma_2 $ is non normal covering of $\Gamma_2.$
\end{center}
{\Pro  Let $\Gamma_n\textcircled{g}\Gamma_r$ with $n > 1$ be the unramified, non normal covering of $\Gamma_r$ then following holds:\\
\begin{enumerate}
\item If $r$ is an even and $av = v$ in $\Gamma_n$ then the alternate action of $a$ and $b^{-1}$ on $y_0 = (v,u_0)$ i.e. 
\begin{eqnarray}
 y_{2p+1} = ay_{2p} , & p = 0,1,\cdots, (2^r-1) \nonumber \\
 y_{2p+2} = b^{-1}y_{2p+1}, & p = 0,1,\cdots, (2^r-2)
\end{eqnarray}
produces spanning path on $v^{th}$ sheet of $\Gamma_n \textcircled{g}\Gamma_r $ which visits every vertex of the $v^{th}$ sheet twice.
\item If $r$ is an odd and $av = v$ in $\Gamma_n$ then the alternate action of $b$ and $a^{-1}$ on $y_0 = (v,u_0)$ i.e. 
\begin{eqnarray}
 y_{2p+1} = by_{2p} , & p = 0,1,\cdots, (2^r-1) \nonumber \\
 y_{2p+2} = a^{-1}y_{2p+1}, & p = 0,1,\cdots, (2^r-2)
\end{eqnarray}
produces spanning path on $v^{th}$ sheet of $\Gamma_n \textcircled{g} \Gamma_r $ which visits every vertex of the $v^{th}$ sheet twice.
\item Let $r$ be any positive number with $av \neq v$ in $\Gamma_n$, then alternate actions of any of the above type goes outside the $v^{th}$ sheet. 
\end{enumerate}} 
\begin{proof}
\begin{enumerate}
\item As $r$ is an even and $av = v $ in $\Gamma_n$, then by definition of generalized replacement product, $v^{th}$ sheet in $\Gamma_n\textcircled{g}\Gamma_r$ contains lift of the edge $e_a$ which connects $(v,u_0)$ to $(v,u_0^a)$ in $\Gamma_n\textcircled{g}\Gamma_r.$  i.e. $v^{th}$ sheet is actually a copy of the graph $\Gamma_r \backslash \{e_b\}.$
So we can write
\begin{eqnarray}
{\bf Rot}_{\Gamma_n\textcircled{g}\Gamma_r}((v,u), a)= ((v,u^{a}), a^{-1}), ~ \forall~ u \in V_{\Gamma_r}  \\
{\bf Rot}_{\Gamma_n\textcircled{g}\Gamma_r}((v,u), b^{-1})= ((v,u^{b^{-1}}), b), ~ \forall u \in V_{\Gamma_r} \backslash \{u_0^{b}\}
\end{eqnarray}
Therefore the simultaneous actions of $a$ and $b^{-1}$ as stated in (10), produces a path which visits the following vertices.
\begin{eqnarray}
y_0 = (v,u_0), \nonumber\\
y_1 = ay_0 = a(v,u_0) = (v,a(u_0)) = (v,u_1), \nonumber\\
y_2 = b^{-1}y_1 = b^{-1}(v,u_1) = (v,b^{-1}u_1) = (v,u_2), \nonumber\\
\ldots
y_{2^{r+1}-1} = ay_{2^{r+1}-2} = (v,u_{2^{r+1}-1})
\end{eqnarray}
One can even prove that $(v,u_{2^{r+1}-1}) = (v,u_0^b)$ i.e. $u_{2^{r+1}-1} = u_0^b$\\
Using induction it has shown as in \cite{DAng}, that action of $b^{-1}a$ has order $2^r$ on $V(\Gamma_r) = X^r.$ Thus $$(b^{-1}a)^{2^r}u_0 = u_0$$ i.e. $$\underbrace{b^{-1}a b^{-1}a \cdots b^{-1}a}_{2^r~~times} u_0 = u_0$$ 
But $$ u_{2^{r+1}-1} = a(b^{-1}a)^{2^{r}-1} u_0 \implies b^{-1}(u_{2^{r+1}-1}) = b^{-1}a(b^-1a)^{2^{r}-1} u_0 = u_0$$
This gives, $$u_{2^{r+1}-1} = u_0^b \implies (v,u_{2^{r+1}-1}) = (v,u_0^b) $$
The right hand side of equation (14), shows that all $y_i's$ are vertices of $v^{th}$ sheet of the graph $\Gamma_n\textcircled{g}\Gamma_r$ which are counted twice. Hence we get the spanning path on $v^{th}$ sheet of $\Gamma_n \textcircled{g} \Gamma_r $ which visits every vertex of the $v^{th}$ sheet twice.
\item As $r$ is an odd number and $av = v $ in $\Gamma_n$, then by definition of generalized replacement product, $v^{th}$ sheet in $\Gamma_n\textcircled{g}\Gamma_r$ contains lift of the edge $e_b$ which connects $(v,u_0)$ to $(v,u_0^b)$ in $\Gamma_n\textcircled{g}\Gamma_r.$  i.e. $v^{th}$ sheet is actually a copy of the graph $\Gamma_r \backslash \{e_a\}.$
So we can write
\begin{eqnarray}
{\bf Rot}_{\Gamma_n\textcircled{g}\Gamma_r}((v,u), b)= ((v,u^{b}), b^{-1}), ~ \forall~ u \in V(\Gamma_r)  \\
{\bf Rot}_{\Gamma_n\textcircled{g}\Gamma_r}((v,u), a^{-1})= ((v,u^{a^{-1}}), a), ~ \forall~ u \in V(\Gamma_r) \backslash \{u_0^{a}\}
\end{eqnarray}
Therefore the simultaneous actions of $b$ and $a^{-1}$ as stated in (11), produces a path which visits the following vertices.
\begin{eqnarray}
y_0 = (v,u_0), \nonumber\\
y_1 = by_0 = b(v,u_0) = (v,b(u_0)) = (v,u_1), \nonumber\\
y_2 = a^{-1}y_2 = a^{-1}(v,u_1) = (v,a^{-1}u_1) = (v,u_2), \nonumber\\
\ldots
y_{2^{r+1}-1} = by_{2^{r+1}-2} = (v,u_{2^{r+1}-1}) = (v,u_0^{a})
\end{eqnarray}
The right hand side of equation (17), shows that all $y_i's$ are vertices of $v^{th}$ sheet of the graph $\Gamma_n\textcircled{g}\Gamma_r$ which are counted twice. Hence we get the spanning path on $v^{th}$ sheet of $\Gamma_n \textcircled{g} \Gamma_r $ which visits every vertex of the $v^{th}$ sheet twice.
\item Let $n > 1$ be positive integer, and $av \neq v \in \Gamma_n.$ The lifts of $e_a$ and $e_b$ go from $v^{th}$ sheet to another sheets. i.e. These lifts connect $(v,u_0)$ to other sheets and hence the $v^{th}$ sheet is actually a copy of the graph $\Gamma_r \backslash \{e_a, e_b\}$, therefore we have
\begin{eqnarray}
{\bf Rot}_{\Gamma_n\textcircled{g}\Gamma_r}((v,u), a^{-1})= ((v,u^{a^{-1}}), a), ~ \forall~ u \in V(\Gamma_r) \backslash \{u_0^{a}\} \\
{\bf Rot}_{\Gamma_n\textcircled{g}\Gamma_r}((v,u), b^{-1})= ((v,u^{b^{-1}}), b), ~ \forall u \in V(\Gamma_r) \backslash \{u_0^{b}\}
\end{eqnarray}
If $r$ is even, the action of $a$ on $(v,u_0)$ gives the vertex $(v^a,u_0^a)$ but $v^a = a(v) \neq v \in \Gamma_n$. Therefore the action of $a$ on $(v,u_0)$ doesn't remain within the $v^{th}$ sheet.
Similarly if $r$ is odd, we can show the action of $b$ on $(v,u_0)$ doesn't remain within the $v^{th}$ sheet. Hence the alternate actions of any of (10) or (11) goes outside the sheet.
\end{enumerate}
\end{proof}
{\Pro If $r > 0$, $n > 1$ and $a(v) \neq v \in \Gamma_n $, then the last $n$ alphabets of the words $a(0^rv)$ and $b(0^rv)$ are distinct from the word $v.$}
\begin{proof}
$$a(0^rv) = \left\{              
                    \begin{array}{ll}           
                     a(0^r)a(v) & \mathrm{if ~r~is~even} \\
                     a(0^r)b(v) & \mathrm{if ~r~is~odd} \\                  
                    \end{array}       
                    \right.$$
But $a(v) \neq v \in \Gamma_n, ~~\forall~ n > 1 $ and $b(v) \neq v , ~~\forall~ v \in \Gamma_n, ~~\forall~ n$
Therefore the last $n$ alphabets of the word $a(0^rv)$ are distinct from the word $v.$
Similarly using action of $b$ on $0^rv$ we can prove that the last $n$ alphabets of the word $b(0^rv)$ are also distinct from the word $v.$
\end{proof}
\section{ZIG ZAG PRODUCT OF GRAPHS} 
\textnormal{The zig zag product of two graphs was introduced by O. Reingold, S. Vadhan and A. Wigderson\cite{Rei}. Taking a product of a large graph with a small graph, the resulting graph inherits (roughly) its size from the large one, its degree from the small one, and its expansion properties from both graphs. The important property of this product is it is a good expander if both large and small graphs are expander. D. D'Angeli, A. Donno and E. Sava-Huss\cite{DAng} have provided the sufficient condition for the zig-zag product of two graphs to be connected.} 
{\Def Let $G_1 = (V_1, E_1)$ be connected $d_1$ regular graph, and let $G_2 = (V_2, E_2)$ be connected $d_2$ regular graph such that $|V_2| = d_1.$ Let ${\bf Rot}_{G_1}$ (resp. ${\bf Rot}_{G_2}$) be the rotation map of $G_1$ (resp. $G_2$). The zig zag product $G_1\textcircled{z}G_2$ is a regular graph of degree $d_2^2$ with vertex set $V_1 \times V_2,$  that we identify with the set $V_1 \times [d_1],$ and whose edges are described by the rotation map 
$$ {\bf Rot}_{G_1\textcircled{z}G_2}((v,k),(i,j)) = ((w,l),(j',l')), $$
for all $ v \in V_1, k \in [d_1], i,j \in [d_2], ~if:$\\
(1) ${\bf Rot}_{G_2}(k,i) = (k',i'),$\\
(2) $ {\bf Rot}_{G_1}(v,k') = (w,l'),$\\
(3) $ {\bf Rot}_{G_2}(l',j) = (l,j'),$\\
where $w \in V_1, l,k',l' \in [d_1]$ and $i',j' \in [d_2]$ and 
$[d_i] = \{1,2,\cdots,d_i\},$ where $ i = 1,2$}\\[8pt]
\textnormal{
D. D'angeli, A. Donno and E. Seva-Huss have shown in Proposition 6.1 of \cite{DAng} that the graph $\Gamma_{n}\textcircled{z}C_4$ is isomorphic to the $double$ $cycle$ $graph~ DC_{2^{n+1}}.$ From the construction of the Schreier graph $\Gamma_n$ of Basilica group see Figure 1, we get the natural labeling of $\Gamma_n$, we label the graph $C_4$ as follows:}\\
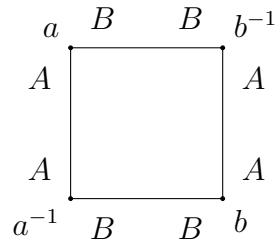
\begin{center}
\begin{tikzpicture}
\draw[] (0,0)node[above right=3pt]  {$B$} -- +(2,0)node[above left=3pt]  {$B$} node[below right=3pt]  {$A$} --+(2,-2)node[above right=3pt]  {$A$} node[below left=3pt]  {$B$} -- +(0,-2) node[below right=3pt]  {$B$} node[above left=3pt]  {$A$}-- +(0,0)node[below left=3pt]  {$A$};
\draw[fill] (0,0) circle [radius=0.025]node[above left]  {$a$};
\draw[fill] (2,0) circle [radius=0.025]node[above right]  {$b^{-1}$};
\draw[fill] (2,-2) circle [radius=0.025]node[below right]  {$b$};
\draw[fill] (0,-2) circle [radius=0.025]node[below left]  {$a^{-1}$};
\end{tikzpicture}
\end{center}
{\Pro  Let $\Gamma_n, n \geq 1$  be the Schreier graph of the action of the Basilica group on $\{0,1\}^n$ and $C_4$ be the $4$ length cycle graph. Let $r \geq 0.$ Then the graph $\Gamma_{n+r}  \textcircled{z} C_4 $ is an unramified graph covering of $\Gamma_{r}  \textcircled{z} C_4 .$ }
\begin{proof}
In the Proposition 2.1, we have shown that the map $\pi : \Gamma_{n+r} \rightarrow \Gamma_r$ is a covering map. Define the map $\Pi : \Gamma_{n+r}  \textcircled{z} C_4  \rightarrow \Gamma_{r}  \textcircled{z} C_4 $ as 
$$ \Pi(uv,s) = (\pi(uv),s) = (u,s), uv \in V({\Gamma_{n+r}}), s \in \{a^{\pm 1}, b^{\pm 1}\}$$
We now show that $\Pi$ is adjacency preserving map.\\
As, $${\bf N}_{\Gamma_{n+r}\textcircled{z} C_4}(uv,a) = \{ ((uv)^{a^{-1}},a^{-1}) , ((uv)^{a^{-1}},b^{-1}) , ((uv)^{b^{-1}},a^{-1}), ((uv)^{b^{-1}},b^{-1})\}$$
So we consider  $(uv,a)$ is adjacent to $((uv)^{b^{-1}},b^{-1})$ and $\Pi(uv,a) = (u,a)$,\\ $\Pi((uv)^{b^{-1}},b^{-1}) = \Pi((u)^{b^{-1}}v',b^{-1}) = ((u)^{b^{-1}},b^{-1}) \in {\bf N}_{\Gamma_{r}\textcircled{z} C_4}(u,a)$\\
Hence, $\Pi(uv,a)$ is adjacent to $\Pi((uv)^{b^{-1}},b^{-1}).$\\
Using a similar argument we can see that other choices of $s$ also gives that  $\Pi$ is adjacency preserving map. \\
Given $(u,s) \in V({\Gamma_{r}\textcircled{z} C_4})$, any vertex $w \in \Pi^{-1}(u,s)$ is of the form $w = (uv,s), $ for some $v \in V({\Gamma_n})$. By definition of the zig zag product of graphs, we see that the vertices in $ V({\Gamma_{r}\textcircled{z} C_4})$ adjacent to $(u,s)$ if $s = a,$ have the form
$${\bf N}_{\Gamma_r \textcircled{z} C_4}(u,a) = \{ (u^{a^{-1}},a^{-1}) , (u^{a^{-1}},b^{-1}) , (u^{b^{-1}},a^{-1}), (u^{b^{-1}},b^{-1})\}$$
and the vertices in $ V({\Gamma_{n+r}\textcircled{z} C_4})$ adjacent to $w = (uv,s) $ if $s = a,$ have the form
$${\bf N}_{\Gamma_{n+r} \textcircled{z} C_4}(uv,a) = \{ ((uv)^{a^{-1}},a^{-1}),
((uv)^{a^{-1}},b^{-1}) , ((uv)^{b^{-1}},a^{-1}), ((uv)^{b^{-1}},b^{-1})\}$$  
Notice that, any $s \in \{a^{\pm 1}, b^{\pm 1}\}, (uv)^s = s(uv).$ The length of the word $s(uv)$ is $n+r.$ By the action of the $s$ on $uv$, the word containing first $r$ alphabets of $s(uv)$ is actually the word $s(u)$ and the remaining $n$ alphabets of $s(uv)$ form some word say $v'$. Thus $\pi((uv)^s) = \pi(s(uv)) = \pi(s(u)v') = s(u) = u^s .$  
$$\implies \Pi({\bf N}_{\Gamma_{n+r} \textcircled{z} C_4}(uv,a)) = {\bf N}_{\Gamma_{r} \textcircled{z} C_4}(u,a)$$
Hence, $\Pi$ is sending neighborhoods of $\Gamma_{n+r} \textcircled{z} C_4$ one to one and onto 
neighborhoods of $\Gamma_{r} \textcircled{z} C_4$ if $s = a.$ Similar argument one can make in the case of $s = b , b^{-1}$ and $ a^{-1}.$
\end{proof}
\textnormal{Notice that if $Y = \Gamma_r \textcircled{z} C_4, $ and $Y' = \Gamma_{1+r} \textcircled{z} C_4,$
\begin{equation} 
{\bf N}_{G}(w,a) = {\bf N}_{G}(w,b), {\bf N}_{G}(w,a^{-1}) = {\bf N}_{G}(w,b^{-1}), ~~for ~~G = Y,Y'\end{equation}
Let $T$ be the spanning subgraph of $Y_1,$ where $V(T) = V({Y})$ and $$E(T) = E({Y}) \backslash \{e_i| i = 1,2,3,4\}$$ where $u_0 = 0^r,$ and, if $r$ is even, then 
$$e_1 = \{(u_0,a^{-1}),(a(b^{-1}a)^{2^r-1}u_0,a)\}, e_2 = \{(u_0,a^{-1}),(a(b^{-1}a)^{2^r-1}u_0,b)\}$$
 $$e_3 = \{(u_0,b^{-1}),(a(b^{-1}a)^{2^r-1}u_0,a)\}, e_4 = \{(u_0,b^{-1}),(a(b^{-1}a)^{2^r-1}u_0,b)\}$$ 
and, if $r$ is odd, then 
$$e_1 = \{(u_0,a^{-1}),(b(a^{-1}b)^{2^r-1}u_0,a)\}, e_2 = \{(u_0,a^{-1}),(b(a^{-1}b)^{2^r-1}u_0,b)\}$$
 $$e_3 = \{(u_0,b^{-1}),(b(a^{-1}b)^{2^r-1}u_0,a)\}, e_4 = 
 \{(u_0,b^{-1}),(b(a^{-1}b)^{2^r-1}u_0,b)\}$$}
{\cor $\Gamma_{1+r}  \textcircled{z} C_4 $ is the two fold normal covering of the graph $\Gamma_{r}  \textcircled{z} C_4 $ with  Galois group $$\displaystyle Gal(\Gamma_{1+r}  \textcircled{z} C_4/\Gamma_r  \textcircled{z} C_4) \cong  \displaystyle \frac{\mathbb{Z}}{2\mathbb{Z}}.$$}
\begin{proof}
The adjacency matrices of graphs  $\Gamma_{k}  \textcircled{z} C_4, k = 1+r,r $ are the circulant matrices, see \cite{DAng}. Therefore there exists the following order of the vertices of  $\Gamma_{1+r}  \textcircled{z} C_4 $  if $r$ is even:\\
\begin{eqnarray}
(0^{1+r},a^{-1});(a(0^{1+r}),a);(b^{-1}a(0^{1+r}),a^{-1});(ab^{-1}a(0^{1+r}),a);\nonumber \\ (b^{-1}ab^{-1}a(0^{1+r}), a^{-1});\cdots(a(b^{-1}a)^{2^{1+r}-1}(0^{1+r}),a);
\end{eqnarray}
\begin{eqnarray}
(0^{1+r},b^{-1});(a(0^{1+r}),b);(b^{-1}a(0^{1+r}),b^{-1});(ab^{-1}a(0^{1+r}),b);\nonumber \\ (b^{-1}ab^{-1}a(0^{1+r}), b^{-1});\cdots (a(b^{-1}a)^{2^{1+r}-1}(0^{1+r}),b).
\end{eqnarray}
In other words we are applying (10) i.e. alternate actions of $a$ and $b^{-1}$ on the word $0^{1+r},$ with the alternation of $a^{-1}, a$ in the second coordinate of the $2^{1+r}$ vertices (the inner cycle of the graph) given in equation (21). \\
Also we are applying same action given in (10) i.e. alternate actions of $a$ and $b^{-1}$ on the word $0^{1+r},$  with the alternation of $b^{-1}, b$ in the second coordinate of $2^{1+r}$ vertices (the outer cycle of the graph) given in equation (22).\\
In the Proposition 2.1, we have also seen that $\Gamma_{n} \textcircled{g} \Gamma_r \simeq \Gamma_{n+r} $ the isomorphism is $f: \Gamma_{n} \textcircled{g} \Gamma_r \rightarrow \Gamma_{n+r}$ defined by  $f(v,u) = uv$ and the isomorphism $f^{-1}: \Gamma_{n+r} \rightarrow \Gamma_{n} \textcircled{g} \Gamma_r$ defined by  $f^{-1}(uv) = (v,u)$, hence we are able to write $ f(0,0^r) = 0^{1+r}$ or $ f^{-1}(0^{1+r}) = (0,0^r)$\\
Thus the vertex $(0^{1+r},a^{-1})$ can be written as $((0,0^r),a^{-1}).$ But the first coordinate of this vertex i.e. $(0,0^r)$ is the $0^r$ vertex of the $0^{th}$ sheet in the graph $\Gamma_1\textcircled{g}\Gamma_r$, and $a(0) = 0 \in \Gamma_1.$ By even case of above Proposition 2.2, we get $$(0,0^r),(0,a(0^r)),(0,b^{-1}a(0^r)),\cdots(0,a(b^{-1}a)^{2^r-1}(0^r))=(0,b(0^r))$$ produces spanning path of the  $0^{th}$ sheet. But again by Proposition 2.1, and $\pi$ is onto map, we have 
\begin{eqnarray}
f(0,0^r) = 0^{1+r} \nonumber  \\
f(0,a(0^r)) =  a(0^{1+r}) \nonumber \\
f(0,b^{-1}a(0^r)) = b^{-1}a(0^{1+r}) \nonumber \\
\cdots f(0,a(b^{-1}a)^{2^r-1}(0^r))= f(0,b(0^r)) = a(b^{-1}a)^{2^r-1}(0^{1+r}) = b(0^{1+r}) \nonumber 
\end{eqnarray}
This gives that the first coordinate of the first $2^r$ vertices of the lists given in equations (21) and (22) contains $0$ in the last position i.e. they are vertices of the $0^{th}$ sheet of $\Gamma_1\textcircled{g}\Gamma_r$. Therefore by definition of zigzag product of graphs, equation (20) and the Proposition 2.2, these vertices forms a (connected) copy of the spanning subgraph $T$ of $Y$. As the first coordinate of these vertices ends with $0$, we denote this copy by $T_0$. Using similar argument, we can show that first coordinate the next $2^r$ vertices given in (21) and (22) ends with $1$ and hence we get the another copy of $T$ which we denote by $T_1$. 
Following are the vertices of $T_0$ and $T_1$.\\
$V({T_0}) = \{ (0^{1+r},a^{-1});(a(0^{1+r}),a);(b^{-1}a(0^{1+r}),a^{-1});(ab^{-1}a(0^{1+r}),a);
\cdots\\
 ~~~~~~~~~~~~~~(a(b^{-1}a)^{2^{r}-1}(0^{1+r}),a);\\  ~~~~~~~~~~~~~~ (0^{1+r},b^{-1});(a(0^{1+r}),b);(b^{-1}a(0^{1+r}),b^{-1});\cdots (a(b^{-1}a)^{2^{r}-1}(0^{1+r}),b) \} \\[5pt]
~~~~~~~~~~~~ = \{(u0,s)| u \in V(\Gamma_r), s \in V(C_4)\} $\\
\vspace{.5cm}\\
$V({T_1}) = \{ ((b^{-1}a)^{2^{r}}0^{1+r},a^{-1});(a(b^{-1}a)^{2^{r}}(0^{1+r}),a);((b^{-1}a)^{2^{r}+1}(0^{1+r}),a^{-1});\cdots\\
 ~~~~~~~~~~~~~~(a(b^{-1}a)^{2^{r+1}-1}(0^{1+r}),a);\\  ~~~~~~~~~~~~~~ ((b^{-1}a)^{2^{r}}0^{1+r},b^{-1});(a(b^{-1}a)^{2^{r}}(0^{1+r}),b);((b^{-1}a)^{2^{r}+1}(0^{1+r}),b^{-1});\\
  ~~~~~~~~~~~~~~\cdots (a(b^{-1}a)^{2^{r+1}-1}(0^{1+r}),b) \} \\[5pt]
~~~~~~~~~~~~ = \{(u1,s)| u \in V(\Gamma_r), s \in V(C_4)\} $\\
and 
$$ V({\widetilde{Y_1}}) = V({T_0}) \cup V({T_1}).$$
Define the non identity Frobenius automorphism map $\varSigma$ as 
$$\varSigma (u1,s) = (\sigma(u1),s) = (u0,s), $$
and
$$\varSigma (u0,s) = (\sigma(u0),s) = (u1,s), ~\forall~ u \in  V({\Gamma_r}),s \in V({C_4}), $$
 where $\sigma \in G(\Gamma_{1+r}|\Gamma_r).$
$$\varSigma(T_0) = T_1, \varSigma(T_1) = T_0 $$ and since $\sigma^2 = id$ 
we have $$ \varSigma^2 = Id.$$
$$(\Pi \circ \varSigma)(u0,s) = \Pi(\varSigma(u0,s)) = \Pi(\sigma(u0),s) = \Pi(u1,s) = (\pi(u1),s) = (u,s)$$
$$(\Pi \circ \varSigma)(u1,s)  = \Pi(\varSigma(u1,s)) = \Pi(\sigma(u1),s) = \Pi(u0,s) = (\pi(u0),s) = (u,s)$$  
and $$\Pi(u0,s) = (u,s) = \Pi(u1,s),$$
$$\implies \Pi \circ \varSigma = \Pi .$$
Thus Galois group $G(Y'|Y) = G(\Gamma_{1+r}\textcircled{z}C_4|\Gamma_r\textcircled{z}C_4)$ contains two elements $\varSigma$ and $\varSigma^2 = Id.$ One can make a similar argument  for the odd $n$. 
\end{proof}
{\cor If $n > 1,$ then $\Gamma_{n+r}\textcircled{z} C_4 $ is the non normal covering of the graph $\Gamma_{r}  \textcircled{z}C_4.$ }\\[8pt]
\textnormal{Proof of this corollary is given at the end of this paper.
To show that $\Gamma_{n+r}\textcircled{z}C_4$ is not a Galois covering of $\Gamma_{r}\textcircled{z}C_4$  we need the following.\\
By definition of zigzag product of graphs,
$$V(\Gamma_{n+r}\textcircled{z}C_4) = \{(uv,s)| uv \in V(\Gamma_{n+r}), s \in V(C_4)\}$$ 
Let $v_0 = 0^n, v_1, v_2, \cdots, v_{2^n - 1}$ be the vertices of the graph $\Gamma_n$.\\
By Proposition 2.1, $\Gamma_{n+r}$ is $2^n$ sheeted cover of $\Gamma_r$, so we can have following partition,
$$ V(\Gamma_{n+r}) = \bigcup_{j = 0}^{2^n -1} \{uv_j | u \in V(\Gamma_r)\}$$
Call $$ S_{v_j} = \{ (uv_j,s) | u \in V(\Gamma_r), s \in V(C_4) \}$$
Therefore we can have
$$V(\Gamma_{n+r}\textcircled{z}C_4) = \bigcup_{j = 0}^{2^n -1} S_{v_j} $$}
\subsection{Restriction Graph}
{\Def Let $G = (V,E)$ be a graph and $ \emptyset \neq V' \subsetneq V$, we define restriction graph $G' = G|_{V'}$ which is a subgraph of $G$ as $G' = (V',E')$\\
$$ E' = \{ e = \{ v_l,v_m\} \in E | v_l, v_m \in V' \}. $$ }
\textnormal{Suppose $$\widetilde{Y} = \Gamma_{n+r}\textcircled{z}C_4, Y = \Gamma_{r}\textcircled{z}C_4 $$
we find the restriction subgraphs $\widetilde{Y}_j$ of $\widetilde{Y}$ by taking $V(\widetilde{Y}_j) = S_{v_j}$\\
One can see, $$V(\widetilde{Y}) = \bigcup_{j = 0}^{2^n -1} V(\widetilde{Y}_j) $$
The next result is about connectedness of restriction subgraphs $\widetilde{Y}_j's$ of $\widetilde{Y}.$}
{\Pro If $\widetilde{Y}_j's$ and $\widetilde{Y}$ as above, then following holds
\begin{enumerate}
\item If $a(v_j) = v_j \in \Gamma_n$ then the restriction subgraph $\widetilde{Y}_j$ of $\widetilde{Y}$ is the connected  $v_j^{th}$ sheet in the graph $\widetilde{Y}.$
\item If $a(v_j) \neq v_j \in \Gamma_n$ then the restriction subgraph $\widetilde{Y}_j$ of $\widetilde{Y}$ is the $v_j^{th}$ sheet in the graph $\widetilde{Y}$ and this sheet is disconnected.
\end{enumerate}} 
\begin{proof}
Recall that $T$ is a spanning subgraph of $Y$, if $r$ is even, we can write 
\begin{eqnarray}
V(T) = \{ (0^r, a^{-1});(a0^r,a);(b^{-1}a0^r,a^{-1}); \cdots;(a(b^{-1}a)^{2^r-1}0^r,a); \nonumber \\
~~~~~~~~~~~(0^r, b^{-1});(a0^r,b);(b^{-1}a0^r,b^{-1}); \cdots;(a(b^{-1}a)^{2^r-1}0^r,b) \} \nonumber
\end{eqnarray}
\begin{enumerate}
\item We write set $T_j, ~~\forall ~~v_j \in \Gamma_n$ as
\begin{eqnarray}
T_j = \{ (0^rv_j, a^{-1});((a(0^r))v_j,a);((b^{-1}a(0^r))v_j,a^{-1}); \cdots;((a(b^{-1}a)^{2^r-1}(0^r))v_j,a); \nonumber \\
~~~~~~~~~~~(0^rv_j, b^{-1});((a(0^r))v_j,b);((b^{-1}a(0^r))v_j,b^{-1}); \cdots;((a(b^{-1}a)^{2^r-1}(0^r))v_j,b) \} \nonumber
\end{eqnarray}
To show that $\widetilde{Y}_j$ is connected it is enough to show that $V(\widetilde{Y}_j) = T_j.$\\
As $(0^rv_j, a^{-1}) = (u_0v_j, a^{-1}) \in S_{v_j} = V(\widetilde{Y}_j)$
and action of $b^{-1}a$ on $0^r \in X^r$ has order $2^r$,
$$ \implies T_j \subset V(\widetilde{Y}_j), ~~\forall ~~j$$
But both the sets $T_j$ and $V(\widetilde{Y}_j)$ contain $2^r$ elements and $T_j \subset V(\widetilde{Y}_j)$. \\ Therefore
$$ T_j = V(\widetilde{Y}_j), ~~\forall ~~j $$
by Proposition 2.1, and equation (20) elements in $T_j$ forms a copy of the spanning subgraph $T$ of $Y$ in $\widetilde{Y}$, and therefore $\widetilde{Y}_j$ is the connected  $v_j^{th}$ sheet in the graph $\widetilde{Y}.$
\item 
By Proposition 2.3, $a(0^rv_j) = a(0^r)a(v_j) \neq a(0^r)v_j$ and  $b(0^rv_j)  \neq b(0^r)v_j.$
and
$$ {\bf N}_{\widetilde{Y}}(0^rv_j, a^{-1}) = \{(a(0^rv_j),a), (a(0^rv_j),b), (b(0^rv_j),a), (b(0^rv_j),b) \} $$
Therefore
$${\bf N}_{\widetilde{Y}}(0^rv_j, a^{-1}) \cap \widetilde{Y}_j = \emptyset, $$ which means $Y_j$ has more than one components in $\widetilde{Y}$, so $\widetilde{Y}_j$ is the $v_j^{th}$ sheet in the graph $\widetilde{Y}$ and this sheet is disconnected. 
\end{enumerate}
\end{proof}
\textnormal{Now we can prove the Corollary 3.2}
\begin{proof}{By contradiction}\\
Let $\widetilde{Y}$ be $2^n$ sheeted normal covering of $Y$. Therefore there exists a Frobenius automorphism $\varSigma \in G(\widetilde{Y}|Y)$ such that $\varSigma$ sends $v_j^{th}$ sheet in $\widetilde{Y}$ to some say $v_k^{th}$ sheet in $\widetilde{Y}$, where $a(v_j) = v_j \in \Gamma_n$ and $a(v_k) \neq v_k \in \Gamma_n.$\\
But by above Proposition 2.4  $v_j^{th}$ sheet is connected and under a Frobenius automorphism it is mapping to disconnected $v_k^{th}$ sheet which is a contradiction.
\end{proof}
\begin{table}[h]
\caption{Vertex labeling for $\Gamma_{1} \textcircled{z} C_4$}
\begin{tabular}{|c|c|c|c|c|c|c|c|c|}
\hline
label & 1 & 2 & 3 & 4 & 5 & 6 & 7 & 8 \\
\hline
\hline
vertex & $(0,a^{-1})$  & $(1,a)$ & $(1,a^{-1})$ & $(0,a)$ &$(0,b^{-1})$  & $(1,b)$ & $(1,b^{-1})$ & $(0,b)$\\
\hline
\end{tabular}
\end{table}
\begin{table}[h]
\caption{Notations for sheets}
\begin{center}
\begin{tabular}{|c|c|}
\hline
Vertex set & Group element\\
\hline
\hline
$\{1,2,\cdots,8\}$ & $\Sigma^2 = Id$\\
$\{1',2',\cdots,8'\}$ & $\Sigma$\\
\hline
\end{tabular}
\end{center}
\end{table}
{\exa The graph $\Gamma_2\textcircled{z}C_4$ is the normal $2$ sheeted cyclic covering of the graph $\Gamma_1\textcircled{z}C_4$ has shown in Figure 5.
We label the vertices using the Table 6. We obtain spanning subgraph of the $\Gamma_{1} \textcircled{z} C_4$ by cutting edges $e_1 = \{1,4\},e_2 = \{1,8\},e_3 = \{5,4\},e_4 = \{5,8\}.$ This gives the dashed line graph $T$ as shown in Figure 6 In the graph $\widetilde{Y}$ there are $16$ vertices of which $8$ ends with the alphabet $0$ and remaining $8$ ends with $1$. By Proposition 2.1, these vertices forms connected sheets which we label in Table 7.\\ 
\begin{figure}
\begin{tikzpicture}[scale=2]
\draw[fill] (45:0.7cm)node[below left] {$1$} circle [radius=0.020] ;
\draw[fill] (135:0.7cm)node[below right ] {$2$}circle [radius=0.020] ;
\draw[fill] (225:0.7cm)node[above right] {$3$}circle [radius=0.020] ;
\draw[fill] (315:0.7cm)node[above left] {$4$}circle [radius=0.020] ;
\draw[fill] (45:1.3cm)node[right] {$5$} circle [radius=0.020] ;
\draw[fill] (135:1.3cm)node[left ] {$6$}circle [radius=0.020] ;
\draw[fill] (225:1.3cm)node[left] {$7$}circle [radius=0.020] ;
\draw[fill] (315:1.3cm)node[right] {$8$}circle [radius=0.020] ;

\draw[thin] (0,0) +(45:0.7cm)   -- +(135:0.7cm) -- +(225:0.7cm) --+(315:0.7cm) --+(45:0.7cm);


\draw[thin] (0,0) +(45:1.3cm)   -- +(135:1.3cm) -- +(225:1.3cm) --+(315:1.3cm) --+(45:1.3cm);


\draw[thin] (0,0) +(45:0.7cm) -- +(135:1.3cm)  -- +(225:0.7cm)--+(315:1.3cm) --+(45:0.7cm);

\draw[thin] (0,0) +(45:1.3cm)  -- +(135:0.7cm) -- +(225:1.3cm) --+(315:0.7cm) --+(45:1.3cm);


\draw[] (270:1.5cm) node[]{${\bf Y = \Gamma_1\textcircled{z}C_4}$}; 

\end{tikzpicture}
\hspace{.8cm}
\begin{tikzpicture}[scale=2]
\draw[fill] (45:0.7cm)node[below left] {$1$} circle [radius=0.020] ;
\draw[fill] (135:0.7cm)node[below right ] {$2$}circle [radius=0.020] ;
\draw[fill] (225:0.7cm)node[above right] {$3$}circle [radius=0.020] ;
\draw[fill] (315:0.7cm)node[above left] {$4$}circle [radius=0.020] ;
\draw[fill] (45:1.3cm)node[right] {$5$} circle [radius=0.020] ;
\draw[fill] (135:1.3cm)node[left ] {$6$}circle [radius=0.020] ;
\draw[fill] (225:1.3cm)node[left] {$7$}circle [radius=0.020] ;
\draw[fill] (315:1.3cm)node[right] {$8$}circle [radius=0.020] ;

\draw[dashed] (0,0) +(45:0.7cm)   -- +(135:0.7cm) -- +(225:0.7cm) --+(315:0.7cm);


\draw[dashed] (0,0) +(45:1.3cm)   -- +(135:1.3cm) -- +(225:1.3cm) --+(315:1.3cm);


\draw[dashed] (0,0) +(45:0.7cm) -- +(135:1.3cm)  -- +(225:0.7cm)--+(315:1.3cm);

\draw[dashed] (0,0) +(45:1.3cm)  -- +(135:0.7cm) -- +(225:1.3cm) --+(315:0.7cm);


\draw[] (270:1.5cm) node[]{${\bf T }$}; 

\end{tikzpicture} \\
\vspace{.25cm}\\
\begin{tikzpicture}[scale=4]

\draw[fill] (0:1cm)node[left] {$1$} circle [radius=0.020] ;

\draw[fill] (315:1cm)node[above left] {$2$}circle [radius=0.020] ;

\draw[fill] (270:1cm)node[above ] {$3$}circle [radius=0.020] ;

\draw[fill] (225:1cm)node[above right] {$4$}circle [radius=0.020] ;

\draw[fill] (180:1cm)node[right] {$1'$}circle [radius=0.020] ;

\draw[fill] (135:1cm)node[below right] {$2'$}circle [radius=0.020] ;

\draw[fill] (90:1cm)node[below ] {$3'$}circle [radius=0.020] ;

\draw[fill] (45:1cm)node[left] {$4'$}circle [radius=0.020] ;

\draw[fill] (0:1.3cm)node[right] {$5$} circle [radius=0.020] ;

\draw[fill] (315:1.3cm)node[below right] {$6$}circle [radius=0.020] ;

\draw[fill] (270:1.3cm)node[below ] {$7$}circle [radius=0.020] ;

\draw[fill] (225:1.3cm)node[below left] {$8$}circle [radius=0.020] ;

\draw[fill] (180:1.3cm)node[left] {$5'$}circle [radius=0.020] ;

\draw[fill] (135:1.3cm)node[above left] {$6'$}circle [radius=0.020] ;

\draw[fill] (90:1.3cm)node[above ] {$7'$}circle [radius=0.020] ;

\draw[fill] (45:1.3cm)node[right] {$8'$}circle [radius=0.020] ;


\draw[dashed] (0,0) +(0:1cm)  -- +(315:1cm)-- +(270:1cm)--+(225:1cm);

\draw[dashed] (0,0) +(180:1cm) --+(135:1cm) -- +(90:1cm) --+(45:1cm);




\draw[dashed] (0,0) +(0:1.3cm)  -- +(315:1.3cm)-- +(270:1.3cm)--+(225:1.3cm);

\draw[dashed] (0,0) +(180:1.3cm) --+(135:1.3cm) -- +(90:1.3cm) --+(45:1.3cm);



\draw[dashed] (0,0) +(0:1.3cm) -- +(315:1cm)-- +(270:1.3cm) --+(225:1cm);
\draw[dashed] (0,0) +(0:1cm) -- +(315:1.3cm)-- +(270:1cm) --+(225:1.3cm);

\draw[dashed] (0,0) +(180:1.3cm) --+(135:1cm) -- +(90:1.3cm) --+(45:1cm);
\draw[dashed] (0,0) +(180:1cm) --+(135:1.3cm) -- +(90:1cm) --+(45:1.3cm);





\draw[thin] (0,0) +(0:1cm) --+(45:1.3cm)--+(0:1.3cm)--+(45:1cm)--+(0:1cm);

\draw[thin] (0,0) +(225:1cm) --+(180:1.3cm)--+(225:1.3cm)--+(180:1cm)--+(225:1cm);



\draw[] (0,0) node[]{${\bf \widetilde{Y} = \Gamma_2\textcircled{z}C_4}$}; 
\end{tikzpicture}
\caption{The graph $\Gamma_2\textcircled{z}C_4$ is  normal $2$ sheeted cyclic covering of the graph $\Gamma_1\textcircled{z}C_4$. The $2$
sheets of $ \widetilde{Y} = \Gamma_2\textcircled{z}C_4$ are copies of the spanning subgraph $T$ of $Y = \Gamma_1\textcircled{z}C_4$ (dashed lines). The continuous lines shown in the graph $\widetilde{Y}$ are lifts of the edges $e_1,e_2,e_3$ and $e_4,$ where $e_1 = \{1,4\}, e_2 = \{ 1,8\}, e_3 = \{5,4\}$ and $e_4 = \{5,8\}.$}
\label{fig3}
\end{figure}
\begin{table}[h]
  \caption{Connections between $1$-sheet  and $\Sigma$-sheet  in $\widetilde{Y} $}
 \begin{center}
 \begin{tabular}{|c|c|}
 \hline
  Vertex & adjacent vertices in $\widetilde{Y} $\\
  \hline
  \hline
  $1$ & $2,6,4',8'$\\
  $2$ & $1,3,5,7$\\
  $3$ & $2,4,6,8$\\
  $4$ & $3,7,1',5'$\\
  $5$ & $2,6,4',8'$\\
  $6$ & $1,3,5,7$\\
  $7$ & $2,4,6,8$\\
  $8$ & $3,7,1',5'$\\
 \hline
 \end{tabular}
 \end{center}
 \end{table}
The representations of the cyclic Galois group $\{1,\Sigma\}$ are the trivial representation $\rho_0$ and the representation defined by $\rho(1) = 1, \rho(\Sigma) = -1.$ So $Q_{\rho} = 2I_8.$
Therefore there are two cases. \\
${\bf Case~~ 1 ~~ The~trivial~representation} ~\rho(1) = {\bf 1}$\\
We use the information given in Table 8 to write the $A(1)$ and $A(\Sigma)$, hence we have
$$ A(1) = \left(
\begin{array}{cccccccc}
 0 & 1 & 0 & 0 & 0 & 1 & 0 & 0 \\
 1 & 0 & 1 & 0 & 1 & 0 & 1 & 0 \\
 0 & 1 & 0 & 1 & 0 & 1 & 0 & 1 \\
 0 & 0 & 1 & 0 & 0 & 0 & 1 & 0 \\
 0 & 1 & 0 & 0 & 0 & 1 & 0 & 0 \\
 1 & 0 & 1 & 0 & 1 & 0 & 1 & 0 \\
 0 & 1 & 0 & 1 & 0 & 1 & 0 & 1 \\
 0 & 0 & 1 & 0 & 0 & 0 & 1 & 0 \\
\end{array}
\right), 
 A(\Sigma) = \left(
\begin{array}{cccccccc}
0 & 0 & 0 & 1 & 0 & 0 & 0 & 1 \\
0 & 0 & 0 & 0 & 0 & 0 & 0 & 0 \\
0 & 0 & 0 & 0 & 0 & 0 & 0 & 0 \\
1 & 0 & 0 & 0 & 1 & 0 & 0 & 0 \\
0 & 0 & 0 & 1 & 0 & 0 & 0 & 1 \\
0 & 0 & 0 & 0 & 0 & 0 & 0 & 0 \\
0 & 0 & 0 & 0 & 0 & 0 & 0 & 0 \\
1 & 0 & 0 & 0 & 1 & 0 & 0 & 0  \\
\end{array}
\right)$$
But  $A_1 = A(1) + A(\Sigma) = A_{\widetilde{Y}}$ ($A_1$ is the adjacency matrix $A_{Y}$ of $Y.$)\\
${\bf  Case~~2~~The ~~representation}~~\rho$\\
Here we find 
$$ A_{\Sigma} = A(1) - A(\Sigma) $$
$$ A_{\Sigma} = \left(
\begin{array}{cccccccc}
 0 & 1 & 0 & -1 & 0 & 1 & 0 & -1 \\
 1 & 0 & 1 & 0 & 1 & 0 & 1 & 0 \\
 0 & 1 & 0 & 1 & 0 & 1 & 0 & 1 \\
 -1 & 0 & 1 & 0 & -1 & 0 & 1 & 0 \\
 0 & 1 & 0 & -1 & 0 & 1 & 0 & -1 \\
 1 & 0 & 1 & 0 & 1 & 0 & 1 & 0 \\
 0 & 1 & 0 & 1 & 0 & 1 & 0 & 1 \\
 -1 & 0 & 1 & 0 & -1 & 0 & 1 & 0 \\
\end{array}
\right)$$
\begin{table}[h]
\caption{Eigenvalues of $A_1$ and $A_{\Sigma}$}
\begin{center}
\begin{tabular}{|c||c|}
\hline
Matrix & Eigenvalues \\
\hline
\hline
$A_1$ & $4,-4,0,0,0,0,0,0 $\\
\hline
$A_{\Sigma}$ & $-2 \sqrt{2},-2 \sqrt{2},2 \sqrt{2},2 \sqrt{2},0,0,0,0$\\
\hline
\end{tabular}
\end{center}
\end{table}
The eigenvalues of $A_1$ and $A_{\Sigma}$ are given in Table 9.\\
\begin{table}[h]
\caption{spectrum of $\widetilde{Y} = \Gamma_2 \textcircled{z} C_4$}
\begin{center}
\begin{tabular}{|c|c|}
\hline
Eigenvalue & Multiplicity\\
\hline
\hline
$\pm 4$ & 1\\
$0$ & 10\\
$\pm 2\sqrt{2}$ & 2 \\
\hline
\end{tabular}
\end{center}
\end{table} 
So we can write the spectrum of $\widetilde{Y} = \Gamma_2 \textcircled{z} C_4.$ See Table 10.\\[5pt]
{\large Reciprocals of $L$ functions for $\Gamma_{2} \textcircled{z} C_4 |\Gamma_{1} \textcircled{z} C_4$.}\\[5pt]
1) For $A_1 $
$$ \zeta_{\Gamma_{1} \textcircled{z} C_4}(t)^{-1} = L(t,A_1,\widetilde{Y}|Y)^{-1} = (1-t^2)^{8} (t-1) (t+1) (3 t-1) (3 t+1) \left(3 t^2+1\right)^6.$$
2) For $A_{\Sigma}$ 
$$ L(t,A_{\Sigma}, \widetilde{Y}|Y)^{-1} = (1-t^2)^{8} \left(3 t^2+1\right)^4 \left(9 t^4-2 t^2+1\right)^2$$
By equation (9) we have, $$\zeta_{\Gamma_{2} \textcircled{z} C_4}(t)^{-1} = L(t,A_1, \widetilde{Y}|Y)^{-1}L(t,A_{\Sigma}, \widetilde{Y}|Y)^{-1}$$ 
$$ = (1-t^2)^{16} (t-1) (t+1) (3 t-1) (3 t+1) \left(3 t^2+1\right)^6 \left(3 t^2+1\right)^4 \left(9 t^4-2 t^2+1\right)^2$$}
\begin{table}[h]
\caption{Notations of sheets }
\begin{center}
\begin{tabular}{|c|c|}
\hline
Vertex set & Sheet index\\
\hline
\hline
$\{1,2,\cdots,8\}$ & $1^{st}$\\
$\{1',2',\cdots,8'\}$ & $2^{nd}$\\
$\{1'',2'',\cdots,8''\}$ & $3^{rd}$\\
$\{1''',2''',\cdots,8'''\}$ & $4^{th}$\\
\hline
\end{tabular}
\end{center}
\end{table}
{\exa Figure 6 contains the example of non normal covering.\\ 
We obtain spanning subgraph of the $\Gamma_{1} \textcircled{z} C_4$ by cutting edges $e_1 = \{1,4\},e_2 = \{1,8\},e_3 = \{5,4\},e_4 = \{5,8\}.$ This gives the dashed line graph $T$ as shown in Figure 6. In the graph $\widetilde{Y}$ there are $16$ vertices of which $8$ ends with the alphabet $0$ and remaining $8$ ends with $1$. The labeling of vertices is given in Table 11. By Proposition 2.1, these vertices forms connected sheets. The labeling of sheets is given in Table 12. 
\begin{table}[h]
\caption{Connections between $2^{nd}$-sheet  in $\widetilde{Y} $ and other sheets }
 \begin{center}
 \begin{tabular}{|c|c|}
 \hline
  Vertex & adjacent vertices in $\widetilde{Y} $\\
  \hline
  \hline
  $1'$ & $4,8,4'',8''$\\
  $2'$ & $1'',3',5'',7'$\\
  $3'$ & $2',4',6',8'$\\
  $4'$ & $3',7',1,5$\\
  $5'$ & $4,8,4'',8''$\\
  $6'$ & $1'',3',5'',7'$\\
  $7'$ & $2',4',6',8'$\\
  $8'$ & $3',7',1,5$\\
\hline
\end{tabular}
\end{center}
\end{table}\\
\begin{figure}
\begin{tikzpicture}[scale=4]

\draw[fill] (0:1cm)node[left] {$1$} circle [radius=0.020] ;
\draw[fill] (337.5:1cm)node[left] {$2$}circle [radius=0.020] ;
\draw[fill] (315:1cm)node[above left] {$3$}circle [radius=0.020] ;
\draw[fill] (292.5:1cm)node[above left] {$4$}circle [radius=0.020] ;
\draw[fill] (270:1cm)node[above ] {$1'$}circle [radius=0.020] ;
\draw[fill] (247.5:1cm)node[above right] {$4''$}circle [radius=0.020] ;
\draw[fill] (225:1cm)node[above right] {$3''$}circle [radius=0.020] ;
\draw[fill] (202.5:1cm)node[right] {$2''$}circle [radius=0.020] ;
\draw[fill] (180:1cm)node[right] {$1'''$}circle [radius=0.020] ;
\draw[fill] (157.5:1cm)node[right] {$2'''$}circle [radius=0.020] ;
\draw[fill] (135:1cm)node[below right] {$3'''$}circle [radius=0.020] ;
\draw[fill] (112.5:1cm)node[below right] {$4'''$}circle [radius=0.020] ;
\draw[fill] (90:1cm)node[below ] {$1''$}circle [radius=0.020] ;
\draw[fill] (67.5:1cm)node[below left] {$2'$}circle [radius=0.020] ;
\draw[fill] (45:1cm)node[left] {$3'$}circle [radius=0.020] ;
\draw[fill] (22.5:1cm)node[left] {$4'$}circle [radius=0.020] ;
\draw[fill] (0:1.3cm)node[right] {$5$} circle [radius=0.020] ;
\draw[fill] (337.5:1.3cm)node[right] {$6$}circle [radius=0.020] ;
\draw[fill] (315:1.3cm)node[below right] {$7$}circle [radius=0.020] ;
\draw[fill] (292.5:1.3cm)node[below right] {$8$}circle [radius=0.020] ;
\draw[fill] (270:1.3cm)node[below ] {$5'$}circle [radius=0.020] ;
\draw[fill] (247.5:1.3cm)node[below left] {$8''$}circle [radius=0.020] ;
\draw[fill] (225:1.3cm)node[below left] {$7''$}circle [radius=0.020] ;
\draw[fill] (202.5:1.3cm)node[left] {$6''$}circle [radius=0.020] ;
\draw[fill] (180:1.3cm)node[left] {$5'''$}circle [radius=0.020] ;
\draw[fill] (157.5:1.3cm)node[left] {$6'''$}circle [radius=0.020] ;
\draw[fill] (135:1.3cm)node[above left] {$7'''$}circle [radius=0.020] ;
\draw[fill] (112.5:1.3cm)node[above left] {$8'''$}circle [radius=0.020] ;
\draw[fill] (90:1.3cm)node[above ] {$5''$}circle [radius=0.020] ;
\draw[fill] (67.5:1.3cm)node[above right] {$6'$}circle [radius=0.020] ;
\draw[fill] (45:1.3cm)node[right] {$7'$}circle [radius=0.020] ;
\draw[fill] (22.5:1.3cm)node[right] {$8'$}circle [radius=0.020] ;

\draw[dashed] (0,0) +(0:1cm) --+(337.5:1cm) -- +(315:1cm)-- +(292.5:1cm);
 \draw[dashed] (0,0) +(247.5:1cm) --+(225:1cm) -- +(202.5:1cm);

\draw[dashed] (0,0) +(180:1cm) --+(157.5:1cm) -- +(135:1cm)--+(112.5:1cm);

\draw[dashed] (0,0) +(67.5:1cm)--+(45:1cm)--+(22.5:1cm) ;


\draw[dashed] (0,0) +(0:1.3cm) --+(337.5:1.3cm) -- +(315:1.3cm)-- +(292.5:1.3cm);
\draw[dashed] (0,0) +(247.5:1.3cm) --+(225:1.3cm) -- +(202.5:1.3cm);

\draw[dashed] (0,0) +(180:1.3cm) --+(157.5:1.3cm) -- +(135:1.3cm)--+(112.5:1.3cm);

\draw[dashed] (0,0)+(67.5:1.3cm)--+(45:1.3cm)--+(22.5:1.3cm) ;

\draw[dashed] (0,0) +(0:1.3cm) --+(337.5:1cm) -- +(315:1.3cm)-- +(292.5:1cm);
\draw[dashed] (0,0) +(247.5:1cm) --+(225:1.3cm) -- +(202.5:1cm);

\draw[dashed] (0,0) +(180:1.3cm) --+(157.5:1cm) --+(135:1.3cm)--+(112.5:1cm);

\draw[dashed] (0,0) +(67.5:1cm)--+(45:1.3cm)--+(22.5:1cm) ;

\draw[dashed] (0,0) +(0:1cm) --+(337.5:1.3cm) -- +(315:1cm)-- +(292.5:1.3cm); 
\draw[dashed] (0,0) +(247.5:1.3cm) --+(225:1cm) -- +(202.5:1.3cm);

\draw[dashed] (0,0) +(180:1cm) --+(157.5:1.3cm) --+(135:1cm)--+(112.5:1.3cm);

\draw[dashed] (0,0) +(67.5:1.3cm)--+(45:1cm)--+(22.5:1.3cm) ;
\draw[thin] (0,0) +(0:1cm) --+(22.5:1.3cm)--+(0:1.3cm)--+(22.5:1cm)--+(0:1cm);

\draw[thin] (0,0) +(202.5:1cm) --+(180:1.3cm)--+(202.5:1.3cm)--+(180:1cm)--+(202.5:1cm);

\draw[thin] (0,0) +(90:1cm) --+(67.5:1cm)--+(90:1.3cm)--+(112.5:1cm)--+(90:1cm) --+(67.5:1.3cm)--+(90:1.3cm)--+(112.5:1.3cm)--+(90:1cm);

\draw[thin] (0,0) +(270:1cm) --+(292.5:1cm)--+(270:1.3cm)--+(247.5:1cm)--+(270:1cm) --+(292.5:1.3cm)--+(270:1.3cm)--+(247.5:1.3cm)--+(270:1cm);

\draw[] (0,0) node[]{${\bf \Gamma_3\textcircled{z}C_4}$}; 
\end{tikzpicture}
\caption{The graph $\Gamma_3\textcircled{z}C_4$ is a non normal $4$ sheeted covering of the graph $\Gamma_1\textcircled{z}C_4$. The $4$ sheets of $\Gamma_3\textcircled{z}C_4$ are copies of the spanning subgraph $T$ of $\Gamma_1\textcircled{z}C_4$ (dashed lines).}
\label{fig}
\end{figure}
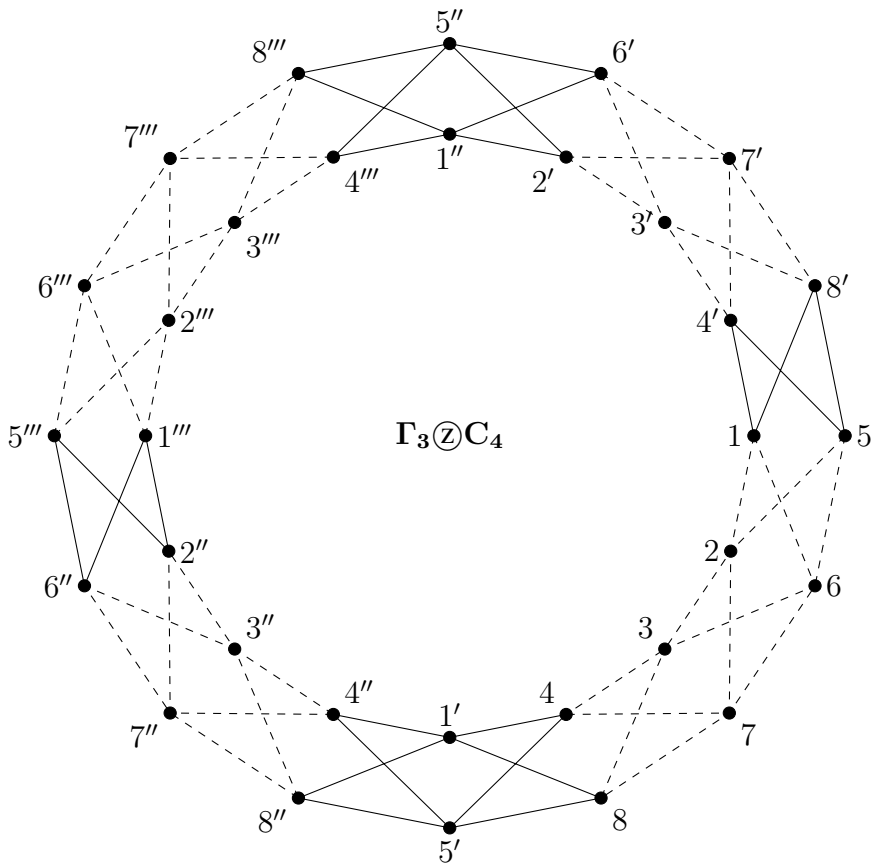
\begin{table}[h]
\caption{Vertex labeling for $\Gamma_{1} \textcircled{z} C_4$}
\begin{tabular}{|c|c|c|c|c|c|c|c|c|}
\hline
label & 1 & 2 & 3 & 4 & 5 & 6 & 7 & 8 \\
\hline
\hline
vertex & $(0,a^{-1})$  & $(1,a)$ & $(1,a^{-1})$ & $(0,a)$ &$(0,b^{-1})$  & $(1,b)$ & $(1,b^{-1})$ & $(0,b)$\\
\hline
\end{tabular}
\end{table}\\
From the Table 13 we can easily observe the following, $$ {\bf N}_{\widetilde{Y}}(1') = \{4,8,4'',8''\} \cap \widetilde{Y_2} = \emptyset$$
By Corollary 3.2, the graph $\Gamma_3 \textcircled{z}C_4$ is non normal $4$ sheeted covering of the graph  $\Gamma_1 \textcircled{z}C_4.$}\\
\textnormal{We propose the following conjecture on the basis of this study. The work presented in this paper may be the special case of this conjecture.}\\
{\conj If $\widetilde{G}$ is normal covering of $G$ with Galois group $Gal(\widetilde{G}|G)$ and $H$ be some other graph, then $\widetilde{G} \textcircled{z}H$ is also normal covering of $G\textcircled{z}H$ with Galois group 
$Gal(\widetilde{G} \textcircled{z}H|G\textcircled{z}H) \simeq Gal(\widetilde{G}|G).$}\\[6pt]
\bibliographystyle{plain}

\end{document}